\documentclass[onecolumn]{IEEEtran}
\usepackage{amsmath,amssymb,epsfig,subfigure,fancybox,balance}
\usepackage{color}
\usepackage{ulem}
\usepackage{graphicx}
\usepackage{enumerate}
\usepackage{cite}
\usepackage{stfloats}

\newcommand{\lbar}{\overline}
\newcommand{\wdh}{\widehat}
\newcommand{\wdt}{\widetilde}

\newcommand{\lf}{\lfloor}
\newcommand{\rf}{\rfloor}
\newcommand{\dl}{\delta}

\def\op{{\mathcal L}}
\def\l{\Big|}
\def\r{\Big|}

\newcommand{\F}{{\mathcal F}}
\newcommand{\e}{\varepsilon}
\newcommand{\rr}{{\Bbb R}}

\newcommand{\cd}{(\cdot)}

\def\para#1{\vskip .2\baselineskip\noindent{\bf #1}}
\def\qed{\strut\hfill $\Box$}

\newtheorem{thm}{Theorem}

\newtheorem{lem}[thm]{Lemma}
\newtheorem{cor}[thm]{Corollary}
\newtheorem{rem}[thm]{Remark}

\newtheorem{exm}[thm]{Example}

\newcommand{\thmref}[1]{Theorem~{\rm \ref{#1}}}
\newcommand{\lemref}[1]{Lemma~{\rm \ref{#1}}}
\newcommand{\corref}[1]{Corollary~{\rm \ref{#1}}}

\newcommand{\remref}[1]{Remark~{\rm \ref{#1}}}
\newcommand{\exmref}[1]{Example~{\rm \ref{#1}}}

\def\al{\alpha}
\def\th{\theta}

\makeatletter
\newcommand{\beq}[1]{\begin{equation} \label{#1}}
\newcommand{\eeq}{\end{equation}}
\newcommand{\bed}{\begin{displaymath}}
\newcommand{\eed}{\end{displaymath}}
\newcommand{\bea}{\bed\begin{array}{rl}}
\newcommand{\eea}{\end{array}\eed}
\newcommand{\ad}{&\!\!\!\disp}
\newcommand{\aad}{&\disp}
\newcommand{\barray}{\begin{array}{ll}}
\newcommand{\earray}{\end{array}}
\def\({\left(}
\def\){\right)}
\def\disp{\displaystyle}

\begin{document}
\date{}
\title{Analyzing Convergence and Rates of Convergence of Particle Swarm Optimization Algorithms Using Stochastic Approximation Methods}
\author{Quan Yuan
and
 George Yin \IEEEmembership{Fellow, IEEE}\thanks{Q. Yuan and G. Yin are with the Department of Mathematics,
 Wayne State University, Detroit, MI 48202, (e-mail: quanyuan@wayne.edu,
 gyin@math.wayne.edu). This research was supported in part by the Army Research Office under grant
W911NF-12-1-0223.}}

\maketitle

\begin{abstract}
Recently, much progress
 has been made on particle swarm optimization (PSO). A number of works
 have been devoted to analyzing the convergence of the underlying algorithms.
Nevertheless, in most cases,
rather simplified hypotheses
are used. For example, it often assumes that the swarm has only one
particle. In addition, more often than not, the variables and the
points of attraction are assumed to remain constant throughout the
optimization process. In reality, such assumptions are often
violated. Moreover, there are no rigorous rates of
convergence results available to date for the particle swarm, to the best of our knowledge. In this paper,
we consider a general form of PSO algorithms, and
analyze asymptotic properties of the algorithms using stochastic
approximation methods. We introduce four coefficients and rewrite
the PSO procedure as a stochastic approximation type iterative
algorithm. Then we analyze its convergence using weak convergence
method. It is proved that a suitably scaled sequence of swarms
converge to the solution of an ordinary differential equation. We
also establish certain stability results. Moreover, convergence
rates are ascertained by using weak convergence method. A centered
and scaled sequence of the estimation errors is shown to have a
diffusion limit.
\end{abstract}


\section{Introduction}\label{int}

\PARstart{R}{ecently}, optimization using particle swarms have received
considerable attention owing to the wide range of applications from
networked systems, multi-agent systems, and autonomous
systems.
Particle swarming refers to
a computational method that optimizes a problem by  trying recursively to improve a candidate solution with respect to a
certain performance measure.
Swarm intelligence from bio-cooperation within groups of
individuals can often provide efficient solutions for certain
optimization problems. When birds are searching food, they exchange
and share information. Each member benefits from all other members
owing to their discovery and experience based on the information
acquired locally.
Then each
 participating member adjusts the next search direction in accordance with the
individual's  best position currently and the information
communicated to this individual by its neighbors.
When food sources scattered unpredictably, advantages of such
collaboration was decisive. Inspired by this, Kennedy and Eberhart
proposed a particle  swarm optimization (PSO) algorithm in 1995
\cite{kennedy19951942}. A PSO procedure is a stochastic optimization
algorithm that mimics the foraging behavior of birds. The search
space of the optimization problem is analogous to the flight space
of birds. Using an abstract setup, each
 bird is
modeled  as a particle (a point in the space of interest). Finding
the optimum is the counterpart of searching for food. A PSO can be
carried out effectively by using an iterative scheme.
The PSO algorithm simulates social behavior
among individuals (particles) ``flying'' through a multidimensional
search space, where each particle represents a point at the
intersection of all search dimensions. The particles evaluate their
positions according to certain fitness functions at each iteration.
The particles share memories of their ``best'' positions locally,
and use the memories to adjust their own velocities and positions.
Motivated by this scenario, a model is proposed to represent the
traditional dynamics of particles.

To put this in a mathematical form, let $F: \rr^D\to \rr$ be the
cost function to be minimized. If we let $M$ denote the size of the
swarm, the current position of particle $i$ is denoted by $X^i$
$(i = 1, 2,\ldots,M)$,
and its current velocity is denoted by $v^i$. Then, the updating
principle can be
expressed as
\beq{eq1}
\barray
  v^{i,d}_{n+1}\ad=v^{i,d}_n+c_1r^{i,d}_{1,n}[\text{Pr}^{i,d}_n-X^{i,d}_n]
  +c_2r^{i,d}_{2,n}[\text{Pg}^{i,d}_n-X^{i,d}_n],\\[0.2cm]
  X^{i,d}_{n+1}\ad=X^{i,d}_n+v^{i,d}_{n+1},
\earray
\eeq
where $d = 1, \ldots, D$; $r^{i,d}_1\sim U(0, 1)$ and $r^{i,d}_2\sim
U(0, 1)$ represent two random variables uniformly distributed in
$[0, 1]$; $c_1$ and $c_2$ represent the acceleration coefficients;
$\text{Pr}^i_n$ represents the best position found by particle $i$
up to ``time'' $n$,
 and $\text{Pg}^i_n$ represents the ``global'' best
position found by particle $i$'s neighborhood $\Pi_i$, i.e.,
\beq{pg-pr}\barray
\ad
\text{Pr}^i_n = \arg\min_{1\leq k\leq n}
F(X^i_k),\\
\ad\text{Pg}^i_n=\text{Pr}^{j^*}_n,\text{  where  }j^*=\mathop{\arg\min}_{j\in\Pi_i}F(\text{Pr}^i_n).\earray\eeq
In
artificial life and social psychology, $v^{i}_n$ in \eqref{eq1} is
the velocity of particle $i$ at time $n$, which provides the
momentum for particles to pass through the search space. The
$c_1r^{i,d}_{1,n}[\text{Pr}^{i,d}_n-X^{i,d}_n]$ is known as the
``cognitive'' component, which represents
 the personal thinking of
each particle. The cognitive component of a particle takes the best
position found so far by this particle as the desired input to make
the particle move toward its own best positions.
$c_2r^{i,d}_{2,n}[\text{Pg}^{i,d}_n-X^{i,d}_n]$ is known as the
``social'' component, which represents the collaborative behavior of
the particles to find the global optimal solution. The social
component always pulls the particles toward the best position found
by its neighbors.

In a nutshell, a PSO algorithm has the following advantages: (1) It
has versatility and does not rely on the problem information; (2) it
has a memory capacity to retain local and global optimal
information; (3) it is easy to implement. Given the versatility and
effectiveness of PSO, it is widely used to solve practical problems
such as artificial neural networks
\cite{juang2004997,messerschmidt2004280}, chemical systems
\cite{costa20031591}, power systems
\cite{alrashidi2009913,abido20011346}, mechanical design
\cite{kovacs2004170}, communications \cite{zhang2004233}, robotics
\cite{li2005554,wu2004147}, economy
\cite{pavlidis2005113,nenortaite2004843}, image processing
\cite{parsopoulos2004402}, bio-informatics
\cite{rasmussen20035,xiao2004895}, medicine \cite{shen2004145}, and
industrial engineering \cite{tasgetiren2004382,liu200718}.
Note that swarms have also been used in many engineering
applications, for example, in collective robotics where
there are teams of robots working together by communicating
over a communication network; see \cite{LiuP04} for a stability analysis and
many related references.

To enable and to enhance further applications, much work  has also
been devoted to improving the PSO algorithms. Because the original
model is similar to a mobile multi-agent system and each parameter
describes a special character of natural swarm behavior, one can
improve the performance of PSO according to the physical meanings of
these parameters
\cite{liu20043751,ratnaweera2004240,poli2005169,parsopouls2002235,zhang20033816}.
The first significant improvement was proposed by Shi and Eberhart
in \cite{shi199869}. They suggested to add a new parameter $w$ as an
``inertia constant'', which results in fast convergence. The
modified equation of \eqref{eq1} is
\begin{equation}\label{eq2}
\begin{aligned}
  v^{i,d}_{n+1}&=wv^{i,d}_n+c_1r^{i,d}_{1,n}[\text{Pr}^{i,d}_n-X^{i,d}_n]+c_2r^{i,d}_{2,n}[\text{Pg}^{i,d}_n-X^{i,d}_n],\\
  X^{i,d}_{n+1}&=X^{i,d}_n+v^{i,d}_{n+1}.
\end{aligned}
\end{equation}
Another significant improvement was due to Clerc and Kennedy
\cite{clerc200258}. They introduced a constriction coefficient
$\chi$ and then proposed to modify \eqref{eq1} as
\begin{equation}\label{eq3}
\begin{aligned}
  v^{i,d}_{n+1}&=\chi(v^{i,d}_n+c_1r^{i,d}_{1,n}[\text{Pr}^{i,d}_n-X^{i,d}_n]+c_2r^{i,d}_{2,n}[\text{Pg}^{i,d}_n-X^{i,d}_n]),\\
  X^{i,d}_{n+1}&=X^{i,d}_n+v^{i,d}_{n+1}.
\end{aligned}
\end{equation}
This constriction coefficient can control the ``explosion'' of the
PSO and ensure the convergence.
Some researchers also considered using ``good'' topologies of
particle connection, in particular adaptive ones (e.g.,
\cite{mendes2004,bratton2007120,miranda2008105}).

There has been much development on mathematical
analysis for the convergence of PSO algorithms as well. Although most
researchers prefer to use discrete system
\cite{trelea2003317,clerc200258,yasuda20031554,brandstaer2002997},
there are some works on
 continuous-time models
\cite{emara20042811,Martinez2011405}. Some recent work such as
\cite{kamisetty2011274,liu20091714,xiao2007396,jiang20078,chen20071271,poli2008}
provides guidelines for selecting PSO parameters
leading to
 convergence, divergence, or oscillation of the swarm's
particles. The aforementioned work also gives rise to several PSO variants.
Nowadays, it is widely recognized that purely
deterministic approach is inadequate in reflecting the exploration
and exploitation aspects brought by stochastic variables.
However, as criticized by Pedersen \cite{peterson2010618}, the
analysis is often oversimplified. For example, the swarm is often assumed to have
only one particle;
stochastic variables
(namely, $r_{1,n}$, $r_{2,n}$) are not used;
the points of attraction,
i.e., the particle's best known position $\text{Pr}$ and the swarm's
best known position $\text{Pg}$,
are normally assumed to remain constant throughout the
optimization process.

In this paper, we study convergence of PSO by using stochastic
approximation methods. In the past,  some authors have considered using stochastic
approximation combined with PSO to enhance the performance or select
parameters (e.g., \cite{kiranyaz200944}). But to the best of our
knowledge, the only paper using stochastic approximation methods to
analyze the dynamics of the PSO so far is by Chen and Li
\cite{chen20071271}. They designed a special PSO procedure and
assumed that
(i) $\text{Pr}^i_n$ and $\text{Pg}^i_n$ are always
within a finite domain;
(ii) with $P^*$ representing
the global optimal positions in the solution
  space,
and $\|P^*\|<\infty$. $\lim_{n\to\infty}\text{Pr}_n\to
  P^*$ and $\lim_{n\to\infty}\text{Pg}_n\to
  P^*$.
Using assumption (i), they proved the  convergence of the algorithm
in the sense of with probability one. With additional assumption
(ii), they showed that the swarm will converge to $P^*$. Despite the
interesting development, their assumptions (i) and (ii) appear to be
rather strong. Moreover, they added some specific terms in the PSO
procedure. So their algorithm is different from the traditional PSOs
\eqref{eq1}-\eqref{eq3}. In this paper, we consider a general form
of PSO algorithms. We introduce four coefficients $\e$, $\kappa_1$,
$\kappa_2$, and $\chi$ and rewrite the PSOs in a stochastic
approximation setup. Then we analyze its convergence using weak
convergence method. We prove that a suitably interpolated sequence
of swarms converge to the solution of an ordinary differential
equation. Moreover, convergence rates are derived by using a
sequence of centered and scaled  estimation errors.

The rest of the paper is arranged as follows. Section
\ref{sec:formulation} presents the setup of our algorithm. Section
\ref{sec:conv} studies the convergence and Section \ref{sec:rate}
analyzes the rate of convergence. Section \ref{sec:exm} proceeds
with several numerical simulation examples to illustrate the
convergence of our algorithms. Finally, Section \ref{sec:rem}
provides a few further remarks.

\section{Formulation}\label{sec:formulation}
First, some descriptions on notation are in order. We use $|\cdot|$
to denote a Euclidean norm. A point $\th$ in a Euclidean space is a
column vector; the $i$th component of $\th$ is denoted by $\th^i$;
$\text{diag}(\th)$ is a diagonal matrix whose diagonal elements are
the elements of $\th$; $I$ denotes the identity matrix of
appropriate dimension; $z'$ denotes the transposition of $z$;
the notation $O(y)$
denotes a function of $y$ satisfying $\sup_{y}|O(y)|/|y|<\infty$,
and $o(y)$ denotes a function of $y$ satisfying $|o(y)|/|y|\to 0$,
as $y\to 0$. In particular, $O(1)$ denotes the boundedness and
$o(1)$ indicates convergence to $0$. Throughout the paper,
for convenience, we use
$K$ to denote a generic positive constant with the convention that
the value of $K$ may be different for different usage.

In this paper, without loss of generality, we assume that each
particle is a one-dimensional scalar. Note that each particle can be
a multi-dimensional vector, which does not introduce essential
difficulties in the analysis; only the notation is a bit more
complex.  We introduce four parameters $\e$, $\kappa_1$, $\kappa_2$,
and $\chi$. Suppose there are $r$ particles, then the PSO algorithm
can be expressed as
\beq{generalform}
\barray
\ad \left[ {\begin{array}{*{20}c}
   {v_{n + 1} }  \\
   {X_{n + 1} }  \\
\end{array} } \right] =  \left[ {\begin{array}{*{20}c}
   {v_n }  \\
   {X_n }  \\
\end{array} } \right]\\
\aad  + \e\left( \left[ {\begin{array}{*{20}c}
   {\kappa _1 I} & { - \chi (c_1 \text{diag}(r_{1,n}) + c_2 \text{diag}(r_{2,n}))}  \\
   {\kappa _2 I} & { - \chi (c_1 \text{diag}(r_{1,n}) + c_2 \text{diag}(r_{2,n}))}  \\
\end{array} } \right]\left[ {\begin{array}{*{20}c}
   {v_n }  \\
   {X_n }  \\
\end{array} } \right]\right.\\
\aad\ \ + \chi\left. \left[ {\begin{array}{*{20}c}
   {c_1 \text{diag}(r_{1,n})}  & c_2 \text{diag}(r_{2,n})   \\
   {c_1 \text{diag}(r_{1,n})}  & c_2 \text{diag}(r_{2,n}) \\
\end{array} } \right]\left[\begin{array}{*{20}c}
{\text{Pr}}(\th_n,\eta_n) \\
{\text{Pg}}(\th_n,\eta_n)\\
\end{array}\right]\right),
\earray
\eeq
\newcounter{tempequationcounter}
\begin{figure*}[!t]
\normalsize
\setcounter{tempequationcounter}{\value{equation}}
\begin{IEEEeqnarray}{rCl}\label{eq6}
\setcounter{equation}{6}
\left[ {\begin{array}{*{20}c}
   {v_{n + 1} }  \\
   {X_{n + 1} }  \\
\end{array} } \right]  =  \left[ {\begin{array}{*{20}c}
   {v_n }  \\
   {X_n }  \\
\end{array} } \right] & + & \e \left\{ \left[ {\begin{array}{*{20}c}
   {\kappa _1 I} & { - 0.5\chi (c_1  + c_2 )I}  \\
   {\kappa _2 I} & { - 0.5\chi (c_1  + c_2 )I}  \\
\end{array} } \right]\left[ {\begin{array}{*{20}c}
   {v_n }  \\
   {X_n }  \\
\end{array} } \right]+ \chi \left[ {\begin{array}{*{20}c}
   {0.5c_1 I} & { 0.5c_2 I}  \\
   {0.5c_1 I} & { 0.5c_2 I}  \\
\end{array} } \right]\left[ {\begin{array}{*{20}c}
   {\text{Pr}(\th_n,\eta_n) }  \\
   {\text{Pg}(\th_n,\eta_n) }  \\
\end{array} } \right]\right.
\nonumber \\
&& +\: \chi \left[ {\begin{array}{*{20}c}
   0 & { - (c_1 \text{diag}(r_{1,n})  + c_2 \text{diag}(r_{2,n})  - 0.5c_1I  - 0.5c_2I )}  \\
   0 & { - (c_1 \text{diag}(r_{1,n})  + c_2 \text{diag}(r_{2,n})  - 0.5c_1I  - 0.5c_2I )}  \\
\end{array} } \right]\left[ {\begin{array}{*{20}c}
   {v_n }  \\
   {X_n }  \\
\end{array} } \right]\nonumber\\
&& +\:
  \left. \chi \left[ {\begin{array}{*{20}c}
   {c_1 \text{diag}(r_{1,n})  - 0.5c_1I} & {c_2 \text{diag}(r_{2,n})  - 0.5c_2I }  \\
   {c_1 \text{diag}(r_{1,n})  - 0.5c_1I} & {c_2 \text{diag}(r_{2,n})  - 0.5c_2I }  \\
\end{array} } \right]\left[ {\begin{array}{*{20}c}
   {\text{Pr}(\th_n,\eta_n) }  \\
   {\text{Pg}(\th_n,\eta_n) }  \\
\end{array} } \right] \right\}. \label{form}
\end{IEEEeqnarray}
\setcounter{equation}
{\value{tempequationcounter}}
\hrulefill \vspace*{4pt}
\end{figure*}
\noindent where
$c_1$ and $c_2$ represent the acceleration coefficients,
$X_n=[X^1_n,\ldots,X^r_n]'\in\rr^r$,
$v_n=[v^1_n,\ldots,v^r_n]'\in\rr^r$,
$\th_n=(X_n,v_n)'$,
 $r_1$, $r_2$ are
$r$-dimensional random vectors in which each component is uniformly
distributed in $(0,1)$, and $\text{Pr}(\th,\eta)$ and
$\text{Pg}(\th,\eta)$ are two non-linear functions depending on
$\theta=(X,v)'$ as well as on a ``noise'' $\eta$, and $\e>0$ is a small
parameter representing the stepsize of the iterations.

\begin{rem}\label{about-cost}
{\rm Note that for a large variety of cases, the structures and the
forms of $\text{Pr}(\th,\eta)$ and $\text{Pg}(\th,\eta)$ are not
known. This is similar to the situation in a stochastic optimization problem in which the
objective function is not known precisely.
Thus, stochastic approximation methods are well suited. As it
is well known that stochastic approximation methods are very useful
for treating optimization problems in which the form of the
objective function is not known precisely, or too complex to
compute. The beauty of such stochastic iteratively defined
procedures is that one need not know the precise form of the
functions.\\
If there is no noise term $\eta_n$, let $\e=0.01$, $\chi=72.9$,
$\kappa_1=-27.1$, and $\kappa_2=72.9$, then \eqref{generalform} is
equivalent to \eqref{eq2} when $w=0.729$ or \eqref{eq3} when
$\chi=0.729$. Thus \eqref{generalform} is a generalization of
\eqref{eq1}-\eqref{eq3}. So a lot of approaches of tuning parameters
(e.g., \cite{beielstein2002,pedersen2010,zhang2005}) could also be
applied.}
\end{rem}

\begin{rem}\label{const-sz}{\rm
In the proposed algorithm, we use a constant stepsize. The stepsize $\e>0$ is
a small parameter. As is well recognized (see \cite{BMP90,kushner2003}), constant stepsize
algorithms have the ability to track slight time variation and is more preferable in many applications.
In the convergence and rate of convergence analysis, we let $\e\to 0$. In the actual computation, $\e$ is just a constant. It need not go to 0.
This is the same as one carries out any computational problem in which the analysis requires the iteration number going to infinity.
However, in the actual computing, one only executes the procedure finitely many steps.
}\end{rem}

In \eqref{generalform}, $r_1$ and $r_2$ are used to reflect the
exploration of particles. Rearranging terms of $\eqref{generalform}$
and considering that $E[c_1\text{diag}(r_{1,n})]=0.5c_1I$ and
$E[c_2\text{diag}(r_{2,n})]=0.5c_2I$, it can be rewritten as
\eqref{form} (on the top of page~\pageref{form}).

Denote
\setcounter{equation}{6}
\beq{eq6a}
\barray
\ad\th_n=[v_n,X_n]'\in\rr^{2r},\\
\ad M = \left[ {\begin{array}{*{20}c}
   {\kappa _1 I} & { - 0.5\chi (c_1  + c_2 )I}  \\
   {\kappa _2 I} & { - 0.5\chi (c_1  + c_2 )I}  \\
\end{array} } \right],\\
\ad P(\th_n,\eta_n)  = \chi{\left[ {\begin{array}{*{20}c}
   {0.5c_1 I} & { 0.5c_2 I}  \\
   {0.5c_1 I} & { 0.5c_2 I}  \\
\end{array} } \right]\left[ {\begin{array}{*{20}c}
   {\text{Pr}(\th_n,\eta_n) }  \\
   {\text{Pg}(\th_n,\eta_n) }  \\
\end{array} } \right]},
\earray
\eeq
and $W(\th_n,r_{1,n},r_{2,n},\eta_n)$ to be the sum of the last
two terms in the curly braces of \eqref{form}. Then \eqref{form} can
be expressed as a stochastic approximation
algorithm
\beq{saform}
\th_{n+1}=\th_n+\e[M\th_n+P(\th_n,\eta_n)+W(\th_n,r_{1,n},r_{2,n},\eta_n)].
\eeq

{One of the challenges in analyzing the convergence of PSO is that
the concrete forms of $\text{Pr}(\th_n,\eta_n)$ and
$\text{Pg}(\th_n,\eta_n)$ are unknown. However, this will not
concern us. As mentioned before, stochastic approximation methods
are known to have advantages in treating such situations.
We shall
use the following assumptions.}

\begin{itemize}
  \item[(A1)] The
  $\text{Pr}(\cdot,\eta)$ and $\text{Pg}(\cdot,\eta)$ are
  continuous
  for each $\eta$.
  For each bounded $\theta$, $E|P(\theta,\eta_n)|^2 <\infty$
  and $E| W(\theta,r_{1,n},r_{2,n},\eta_n)|^2 <\infty$.
  There exist continuous functions
  $\lbar{\text{Pr}}(\th)$ and $\lbar{\text{Pg}}(\th)$ such
  that
   \beq{eq24-0}
  \barray
  \ad\frac{1}{n}\sum^{n+m-1}_{j=m}E_m\text{Pr}(\th,\eta_j)\to\lbar{\text{Pr}}(\th) \ \hbox{ in probability,}\\
  \ad\frac{1}{n}\sum^{n+m-1}_{j=m}E_m\text{Pg}(\th,\eta_j)\to\lbar{\text{Pg}}(\th) \ \hbox{ in probability,}\\
 \earray
  \eeq
where
$E_m$ denotes the conditional expectation on the $\sigma$-algebra
$\F_m=\{\th_0, r_{i,j}, i=1,2, \eta_j: j< m\}$. Moreover, for each
$\theta$ in a bounded set,
  \beq{eq24}
  \barray
  \disp\sum^{\infty}_{j=n}|E_n\text{Pr}(\th,\eta_j)-\lbar{\text{Pr}}(\th)|\ad<\infty,\\
  \disp\sum^{\infty}_{j=n}|E_n\text{Pg}(\th,\eta_j)-\lbar{\text{Pr}}(\th)|\ad<\infty.
  \earray
  \eeq

  \item[(A2)] Define $$\lbar{P}(\theta)= \chi{\left[ {\begin{array}{*{20}c}
   {0.5c_1 I} & { - 0.5c_2 I}  \\
   {0.5c_1 I} & { - 0.5c_2 I}  \\
\end{array} } \right]\left[ {\begin{array}{*{20}c}
   {\lbar{\text{Pr}}(\th) }  \\
   {\lbar{\text{Pg}}(\th) }  \\
\end{array} } \right]}.$$ The ordinary differential equation
  \beq{ODE}
  \frac{d\th(t)}{dt}=M\th(t)+\lbar{P}(\th(t))
  \eeq
  has a unique solution for each initial
  condition $\th(0)=(\th^1_0,\ldots,\th^{2r}_0)'$.

\item[(A3)]  For $i=1,2$, $\{r_{i,n}\}$ and $\{\eta_n\}$ are
 mutually independent; $\{r_{i,n} \}$ are i.i.d. sequences of random variables
with each component
 being uniformly distributed in $(0,1)$.

\end{itemize}

\begin{rem}\label{rem:about-a2} {\rm
Condition (A1) is satisfied by a large class of functions and random variables.
The continuity is assumed for convenience. In fact, only weak continuity is needed so we can in fact
deal with indicator type of functions whose expectations are continuous.

In fact, \eqref{eq24-0} mainly requires
that $\{\text{Pr}(\th,\eta_n)\}$ is a sequence that satisfies a law of large number type of condition, although it is weaker than the
usual weak law of large numbers. Condition
\eqref{eq24} is modeled by the mixing type condition.
For instance, we may assume that
for each bounded random vector $\theta$ and each $T<\infty$,
either
\bea \ad \lim_{j\to \infty, \Delta\to 0} E \sup_{|Y|\le \Delta} |\text{Pr}(\th+Y ,\eta_j)-\text{Pr}(\th,\eta_j)|= 0, \ \hbox{ or}\\
\ad \lim_{n\to \infty, \Delta\to 0} {1\over n} \sum^{m+n-1}_{j=m} E \sup_{|Y|\le \Delta} |\text{Pr}(\th+Y ,\eta_j)-\text{Pr}(\th,\eta_j)|= 0.\eea
Apparently, the second alternative is even weaker. With either of this assumption,
all of the subsequent development follows, but the argument is more complex.
Under the above condition, one can treat discontinuity involving sign function or indicator function among others.
For the corresponding stochastic approximation algorithms, see \cite[p. 100]{Kushner84};
the setup in \cite{kushner2003} is even more general, which allows in addition to the discontinuity,
the functions involved to be time dependent.
Inserting the conditional expectation is much weaker than without.
For example,  if $\{\eta_n\}$ is a sequence of i.i.d. random variables with distribution function
$F_\eta$, then for each $\th$, $\lbar{\text{Pr}}(\th) =
E \text{Pr}(\th, \eta_1) =
\int \text{Pr}(\th, \zeta) F_\eta(d\zeta)$, so \eqref{eq24-0} is easily verified.
Likewise, if $\{\eta_n\}$ is a  martingale difference sequence, the condition is also satisfied. Next, if $\{\eta_n \}$ is a moving average sequence driven by a martingale difference sequence, \eqref{eq24-0} is also satisfied. In addition, if $\{\eta_n\}$ is a mixing sequence \cite[p.166]{Billingsley} with the mixing rate decreasing to 0, the condition is also satisfied.
Note that in a mixing sequence, there can be infinite correlations and the remote past and distant future are only asymptotically uncorrelated.

In the simplest additive noise case, i.e.,
 $\text{Pr}(\th,\eta) =\text{Pr}(\th)+ \eta$, then the condition is mainly on the
noise sequence $\{\eta_n\}$.
Condition \eqref{eq24} is modeled after the so-called mixing inequality; see \cite[p.82]{Kushner84} and references therein.
Suppose that $\{\text{Pr}(\th,\eta_n)\}$ is a stationary mixing sequence with mean $\lbar{\text{Pr}}(\th)$ and mixing rate
$\phi_n$ such that $\sum_n \phi_n^{1/2} <\infty$, then \eqref{eq24} is satisfied.

{With these assumptions, we proceed to analyze the convergence and
rates of convergence  of PSO algorithms with general form
\eqref{saform}. The scheme is a constant-step-size stochastic
approximation algorithm with step size $\e$. Our interest lies in
obtaining convergence and rates of convergence as $\e\to0$. We
emphasize that in the actual computation, it is not necessary to
modify it as the generalized PSO form \eqref{saform}. This
generalized PSO form is simply a convenient form that allows us to
analyze the algorithm by using methods of stochastic
approximation.}}\end{rem}

\section{Convergence}\label{sec:conv}

This section is devoted to obtaining asymptotic properties of
 algorithm \eqref{saform}. In relation to PSO the word
``convergence''
typically means one of two things, although it is often
not clarified which definition is meant and sometimes they are
mistakenly thought to be identical.
(i) Convergence may refer to the swarm's best known position $\rm Pg$ approaching (converging to) the optimum of
  the problem, regardless of how the swarm behaves.
(ii) Convergence may refer to a swarm collapse in which all particles have converged to a point in the search space,
  which may or may not be the optimum.
 Since the convergence may
rely on structure of the cost function if we use the first
definition of convergence, we use the second one as the definition
of convergence in this study. Our first result concerns the property
of the algorithm as $\e\to 0$ through an appropriate continuous-time
interpolation. We define
$$\th^\e(t)=\th_n\ \hbox{ for } \ t\in[\e n,\e n+\e).$$ Then
$\th^\e\cd\in D([0,T]:\rr^{2r})$, which is the space of functions
that are defined on $[0,T]$ taking values in $\rr^{2r}$, and that
 are right continuous and have left limits endowed
with the Skorohod topology \cite[Chapter 7]{kushner2003}.

\begin{thm}\label{conv}
Under {\rm(A1)-(A3)},
$\th^\e\cd$ is tight in $D([0,T]: \rr^{2r})$. Moreover, as $\e\to
0$, $\th^\e\cd$ converges weakly to $\th\cd$, which is a solution of
\eqref{ODE}.
\end{thm}

\begin{rem}{\rm
An equivalent way of stating the ODE limit \eqref{ODE}
is to consider its associated martingale problem \cite[pp. 15-16]{Kushner84}. Consider the
differential operator associated with $\th\cd$
\[\op f(\th)=(\nabla f(\th))'(M\th+\lbar{P}(\th)),\]
and
define
\[\wdt M_f(t)=f(\th(t))-f(\th(0))-\int^t_0\op f(\th(s))ds.\]
If $\wdt M_f\cd$ is a martingale for each $f\cd\in C^1_0$ ($C^1$
function with compact support), then $\th\cd$ is said to solve a
martingale problem with operator $\op$. Thus, an equivalent way to
state the theorem is to prove that $\th^\e\cd$ converges weakly to
$\th\cd$, which is a solution of the martingale problem with
operator $\op$. }\end{rem}

\para{Proof of \thmref{conv}.}
To prove the tightness in $D([0,T]: \rr^{2r})$, we first need to show
\beq{eq30}
\lim_{K\to\infty}\limsup_{\e\to0} P\{\sup_{t\leq T}|\th^\e(t)|\geq K\}=0
\eeq
To avoid verifying \eqref{eq30}, we define a process $\th^{\e,N}\cd$
satisfies $\th^{\e,N}(t)=\th^\e(t)$ up until the first exit from
$S_N=\{ x\in \rr^{2r}: | x| \le N\}$ and satisfies \eqref{eq30}, the
$\th^{\e,N}\cd$ is said to be an \emph{$N$-truncation} of
$\th^\e\cd$. Introduce a truncation function $q^N\cd$ that is
smooth
and that satisfies $q^N(\th)=1$ for
$|\th|\leq N$, $q^N(\th)=0$ for $|\th|\geq N+1$. Then
the discrete system \eqref{saform} is defined as
\beq{saform-tru}
\barray
\th^N_{n+1}\ad=\th^N_n+\e[M\th^N_n+P(\th^N_n,\eta_n)\\
\aad\qquad\qquad\quad+W(\th^N_n,r_{1,n},r_{2,n},\eta_n)]q^N(\th^N_n),
\earray
\eeq
using the $N$-truncation. Moreover, the $N$-truncated ODE and the
operator $\op^N$ of the associated martingale problem can be
defined as
\beq{ODE-tru}
{d \theta^N(t) \over dt}
=[M\th^N(t)+\lbar P(\th^N(t))]q^N(\th(t)),
\eeq
and
\beq{op-tru}
\op^N f(\th)=(\nabla f(\th))'[M\th+\lbar P(\th)]q^N(\th),
\eeq
respectively.

To prove the theorem,  we proceed to verify the following claims:
(a) for each $N$, $\{\th^{\e,N}\cd\}$ is tight. By virtue of the
Prokhorov theorem \cite[p.229]{kushner2003}, we can extract a weakly
convergent subsequence. For notational simplicity, we still denote the
subsequence by $\{\th^{\e,N}\cd\}$ with limit denoted by $\th^N\cd$.

(b) $\th^N\cd$ is a solution of the martingale problem with operator
$\op^N$.
 Using the uniqueness of the limit, passing to the limit as
$N\to \infty$, and
 by the corollary in \cite[p.44]{Kushner84}, $\{\th^\e\cd\}$
converges weakly to $\th\cd$.

Now we start to prove claims (a) and
(b).

(a) Tightness. For any $\dl>0$, let $t>0$ and $s>0$ such that
$s\leq\dl$, and $t$,
$t+\dl\in[0,T]$. Note that
\bea
\ad
\th^{\e,N}(t+s)-\th^{\e,N}(t)\\
\aad =\e\sum^{(t+s)/\e-1}_{k=t/\e}(M\th^N_k+P(\th^N_k,\eta_k)
\\
\aad\qquad\qquad\qquad+W(\th^N_k,r_{1,k},r_{2,k},\eta_k))q^N(\th^N_k).
\eea
In the above and hereafter, we use the conventions that $t/\e$ and
$(t+s)/\e$ denote the corresponding integer parts $\lf t/\e \rf$ and
$\lf (t+s)/\e \rf$, respectively. For notational simplicity, in what
follows, we will not use the floor function notation unless it is
necessary.

Using the Cauchy-Schwarz inequality,
\beq{eq32}
\barray
\ad\e^2E^\e_t\left|\sum^{(t+s)/\e-1}_{k=t/\e}M\th^N_kq^N(\th^N_k)\right|^2\\
\aad\quad \leq
\e Ks\sum^{(t+s)/\e-1}_{k=t/\e}E^\e_t\left|\th^N_kq^N(\th^N_k)\right|^2.
\earray
\eeq
where $E^\e_t$ denotes the expectation
conditioned on the $\sigma$-algebra $\F^\e_t$. Likewise,
\beq{eq33}
\e^2E^\e_t\left|\sum^{(t+s)/\e-1}_{k=t/\e}
W(\th^N_k,r_{1,k},r_{2,k},\eta_k)q^N(\th^N_k)\right|^2\leq Ks^2,
\eeq
and
\beq{eq34}
\e^2E^\e_t\left|\sum^{(t+s)/\e-1}_{k=t/\e}P(\th^N_k,\eta_k)q^N(\th^N_k)\right|^2\leq Ks^2.
\eeq
So we have
\beq{eq312}
\barray
\ad E^\e_t\l \th^{\e,N}(t+s)-\th^{\e,N}(t)\r^2\\
\aad \leq K\e s\sum^{(t+s)/\e-1}_{k=t/\e}\sup_{t/\e\leq k\leq
(t+s)/\e-1} E^\e_t|\th^N_kq^N(\th^N_k)|^2+Ks^2.
\earray
\eeq
As a result,  there is a $\varsigma^\e(\delta)$ such that
$$E^\e_t|\th^{\e,N}(t+s)-\th^{\e,N}(t)|^2 \le E^\e_t \varsigma^\e(\delta) \ \hbox{ for all } \ 0\le s\le \delta,$$
and that
$\lim_{\dl\to 0}\limsup_{\e\to
0}E \varsigma^\e(\delta)=0.$ The tightness
of $\{\th^{\e,N}\cd\}$ then follows from \cite[p.47]{Kushner84}.

(b) Characterization of the limit.
To characterize the limit process, we need to work with a
continuously differentiable function with compact support $f\cd$.
Choose $m_\e$ so that $m_\e\to\infty$ as $\e\to 0$ but $\dl_\e=\e m_\e\to 0$.
Using the recursion \eqref{saform-tru},
\beq{eq39}
\barray
\ad\!\!\!\! f(\th^{\e,N}(t+s))-f(\th^{\e,N}(t))\\
\ad\! =\sum^{(t+s)/\dl_\e}_{l=t/\dl_\e}[f(\th^N_{lm_\e+m_\e})-f(\th^N_{lm_\e})]\\
\ad\! =\e\!\!\sum^{(t+s)/\dl_\e}_{l=t/\dl_\e}(\nabla
f(\th^N_{lm_\e}))'
\sum^{lm_\e+m_\e-1}_{k=lm_\e}[M\th^N_k+\lbar{P}(\th^N_k)]q^N(\th^N_k)\\
\ad +\e\!\!\!\sum^{(t+s)/\dl_\e}_{l=t/\dl_\e}\!\!\!(\nabla
f(\th^N_{lm_\e}))'
\!\!\sum^{lm_\e+m_\e-1}_{k=lm_\e}\!\![P(\th^N_k,\eta_k)-\lbar{P}(\th^N_k)]q^N(\th^N_k)\\
\ad +\e\!\!\!\sum^{(t+s)/\dl_\e}_{l=t/\dl_\e}\!\!\!(\nabla
f(\th^N_{lm_\e}))'
\sum^{lm_\e+m_\e-1}_{k=lm_\e}\!\!\!W(\th^N_k,r_{1,k},r_{2,k},\eta_k)q^N(\th^N_k)\\
\ad +\e\!\!\sum^{(t+s)/\dl_\e}_{l=t/\dl_\e}\Bigg\{(\nabla
f(\th^{N+}_{lm_\e})-\nabla
f(\th^{N}_{lm_\e}))'
\\
\aad\quad\quad\times
\sum^{lm_\e+m_\e-1}_{k=lm_\e}[M\th^N_k+P(\th^N_k,\eta_k) \\
\aad \qquad\qquad \qquad +W(\th^N_k,r_{1,k},r_{2,k},\eta_k)]q^N(\th^N_k)\Bigg\},
\earray
\eeq
where $\th^{N+}_{lm_\e}$ is a point on the line segment joining
$\th^{N}_{lm_\e}$ and $\th^{N}_{lm_\e+m_\e}$.

Our focus here is to characterize the limit. By the Skorohod
representation \cite[p.230]{kushner2003}, with a slight abuse of
notation, we may assume that $\th^{\e,N}\cd$ converges to $\th^N\cd$
with probability one and the convergence is uniform on any bounded
time interval. To show that $\{\th^{\e,N}\cd\}$ is a solution of the
martingale problem with operator $\op^N$, it suffices to show that
for any $f\cd \in C^1_0$, the class of functions that are
continuously differentiable with compact support,
$$\wdt M^N_f(t)=f(\th^N(t))-f(\th^N(0))-\int^t_0\op^N f(\th^N(u))du$$ is a
martingale. To verify the martingale property, we need only show
that for any bounded and continuous function $h\cd$, any positive
integer $\kappa$, any $t,\,s
> 0$, and $t_i\leq t$ with $i\leq\kappa$,
\beq{eq310}
\barray
\ad Eh(\th^N(t_i):i\leq\kappa)[\wdt M^N_f(t+s)-\wdt M^N_f(t)]\\
\aad =
Eh(\th^N(t_i):i\leq\kappa)\\
\aad\quad\times[f(\th^N(t+s))-f(\th^N(t))-\int^{t+s}_t\op^N
f(\th^N(u))du]\\
\aad =0.
\earray
\eeq
To verify \eqref{eq310}, we begin with the process indexed by $\e$.
For notational simplicity, denote
\beq{eq311}
\wdt h=h(\th^N(t_i):i\leq\kappa),\ \ \wdt h^\e=h(\th^{\e,N}(t_i):i\leq\kappa).
\eeq
Then the weak convergence and the Skorohod representation together
with the boundedness and the continuity of $f\cd$ and $h\cd$ yield
that as $\e\to0$,
\bea \ad E\wdt h^\e[f(\th^{\e,N}(t+s))-f(\th^{\e,N}(t))\\
\aad \to E\wdt h[f(\th^N(t+s))-f(\th^N(t))].\eea

\def\bgl{\Big[}
\def\bgr{\Big]}
For the last term of \eqref{eq39}, as $\e\to0$, since $f\cd\in C^1_0$,
\beq{eq312-0}
\barray
\ad\  E\wdt h^\e \e\sum^{(t+s)/\dl_\e}_{l=t/\dl_\e}\Bigg\{(\nabla
f(\th^{N+}_{lm_\e})-\nabla
f(\th^{N}_{lm_\e}))'\\
\aad\qquad\qquad\ \times\sum^{lm_\e+m_\e-1}_{k=lm_\e}[M\th^N_k+P(\th^N_k,\eta_k)\\
\aad \qquad\qquad\       +W(\th^N_k,r_{1,k},r_{2,k},\eta_k)]q^N(\th^N_k)\Bigg\}\\
\aad\qquad =O(\e)\to0.
\earray
\eeq
For the next to the last term of \eqref{eq39},
\beq{eq312-aa}
\barray
\ad\!\!\! \lim_{\e\to0} E\wdt h^\e\Big[\e\sum^{(t+s)/
\dl_\e}_{l=t/\dl_\e}(\nabla f(\th^N_{lm_\e}))'\\
\aad\quad\quad\times\sum^{lm_\e+m_\e-1}_{k=lm_\e}W(\th^N_k,r_{1,k},
r_{2,k},\eta_k)q^N(\th^N_k)\Big]\\
\aad \ =\lim_{\e\to0} E\wdt h^\e\Big[\sum^{(t+s)/\dl_\e}_{l=t/\dl
_\e}(\nabla f(\th^N_{lm_\e}))'\\
\aad\quad\quad\times\frac{\dl_\e}{m_\e}
\sum^{lm_\e+m_\e-1}_{k=lm_\e}
E_{lm_\e}W(\th^N_k,r_{1,k},r_{2,k},\eta_k)q^N(\th^N_k)\Big].
\earray
\eeq
Using (A1) and (A3),
\[\frac{1}{m_\e}\sum^{l m_\e + m_\e -1}_{j=lm_\e}E_{lm_\e}W(\th^N_{lm_\e},
r_{1,j},r_{2,j},\eta_j)q^N(\th^N_{lm_\e})\to0\]
in probability,
we obtain that
\beq{eq313}
\barray
\ad E\wdt h^\e\Big[\e\sum^{(t+s)/\dl_\e}_{l=t/\dl_\e}(\nabla
f(\th^N_{lm_\e}))'\\
\aad\quad\quad\times\sum^{lm_\e+m_\e-1}_{k=lm_\e}W(\th^N_k,
r_{1,k},r_{2,k},\eta_k)q^N(\th^N_k)\Big]\to0.
\earray
\eeq
Using (A1), we obtain
\beq{eq314}
\barray
\ad E\wdt h^\e\bgl\e\sum^{(t+s)/\dl_\e} _{l=t/\dl_\e}(\nabla
f(\th^N_{lm_\e}))'\\
\aad\quad\quad\times\sum^{lm_\e+m_\e-1}
_{k=lm_\e}(P(\th^N_k,\eta_k)-\lbar{P}(\th^N_k))q^N(\th^N_k)\bgr
\to0.
\earray
\eeq
Next, we consider the first term. We have
\beq{eq315}
\barray
\ad \lim_{\e\to0}E\wdt
h^\e\bgl\e\sum^{(t+s)/\dl_\e}_{l=t/\dl_\e}(\nabla
f(\th^N_{lm_\e}))'\\
\aad\quad\quad\times\sum^{lm_\e+m_\e-1}_{k=lm_\e}
(M\th^N_k+\lbar{P}(\th^N_k))q^N(\th^N_k)\bgr\\
\aad = \lim_{\e\to0}E\wdt
h^\e\bgl\e\sum^{(t+s)/\dl_\e}_{l=t/\dl_\e}(\nabla
f(\th^N_{lm_\e}))'\\
\aad\quad\quad\times\sum^{lm_\e+m_\e-1}_{k=lm_\e}
(M\th^N_{lm_\e}+\lbar{P}(\th^N_{lm_\e}))q^N(\th^N_{lm_\e})\bgr.
\earray
\eeq
Thus, to get the desired limit, we need only examine the last two lines
above. Let $\e lm_\e\to u$ as $\e\to0$. Then for all $k$ satisfying
$lm_\e\leq k\leq lm_\e+m_\e-1$, $\e k\to u$ since
$\dl_\e\to0$. As a result,
\beq{eq316}
\barray
\ad \lim_{\e\to0}E\wdt h^\e\bgl
\e\sum^{(t+s)/\dl_\e}_{l=t/\dl_\e}(\nabla
f(\th^N_{lm_\e}))'\\
\aad\quad\quad\times\sum^{lm_\e+m_\e-1}_{k=lm_\e}
(M\th^N_{lm_\e}+\lbar{P}(\th^N_{lm_\e}))q^N(\th^N_{lm_\e})\bgr\\
\aad = E\wdt h\bgl\int^{t+s}_t(\nabla
f(\th^N(u)))'(M(\th^N(u))\\
\aad\qquad\quad+\lbar{P}(\th(u)))q^N(\th(u))du\bgr.
\earray
\eeq
The desired result then follows.\qed

To proceed, consider \eqref{ODE}. For simplicity,
suppose that there is a unique stationary point $\theta^*$.
 Denote
$\lbar{\text{Pr}}(\th^*)=\text{Pr}^*$ and
$\lbar{\text{Pg}}(\th^*)=\text{Pg}^*$.
By the inversion formula of partitioned matrix \cite{tian20091294},
solving $M\theta^*+ \lbar P(\theta^*) =0$
yields that the equilibrium point of the ODE satisfies
\beq{eq317}
\barray
\th^* \ad = \left[ {\begin{array}{*{20}c}
   {\kappa _1 I} & { - 0.5\chi (c_1  + c_2 )I}  \\
   {\kappa _2 I} & { - 0.5\chi (c_1  + c_2 )I}  \\
 \end{array} } \right]^{ - 1}\\
 \aad\quad\quad \times
 \left[ {\begin{array}{*{20}c}
   { - 0.5\chi (c_1 \text{Pr}^*  + c_2 \text{Pg}^* )}  \\
   { - 0.5\chi (c_1 \text{Pr}^*  + c_2 \text{Pg}^* )}  \\
\end{array} } \right]\\
\ad= \left[ {\begin{array}{*{20}c}
   0  \\
   {\frac{{c_1 \Pr ^*  + c_2 \text{Pg}^* }}
{{c_1  + c_2 }}}  \\
\end{array} } \right].
\earray
\eeq

\begin{cor}\label{conv-to-t*}
Suppose that the stationary point $\theta^*$ is asymptotically stable in the sense of Lyapunov and that
$\{\th_n\}$ is tight.
Then for any $t_\e\to \infty$ as $\e\to 0$,
$\theta^\e(\cdot + t_\e)$ converges weakly to $\theta^*$.
\end{cor}

\para{Proof.} Define $\wdt \th^\e\cd = \th^\e(\cdot +t_\e)$.
Let $T> 0$ and
consider the pair $\{\wdt\theta^\e(  \cdot) ,\wdt \th^\e (\cdot-T )\}$.
Using the same argument as in the proof of \thmref{conv},
$\{\wdt\theta^\e(  \cdot) , \wdt\th^\e (\cdot - T )\}$ is tight.
Select a
convergent subsequence with limit
denoted by $(\theta\cd, \theta_T\cd)$.
Then $\th(0)= \th_T(T)$. The value of
$\th_T(0)$ is not known, but all  such $\th_T(0)$, over all $T$ and convergent subsequences,
belong to a tight set. This together with the stability
 and \thmref{conv} implies that for any $\Delta > 0$, there is a $T_\Delta$ such that
for $T > T_\Delta$, $ P(\th_T(T) \in U_\Delta (\th^*) ) > 1- \Delta$,
where $U_\Delta(\th^*)$ is a $\Delta$-neighborhood of
$\th^*$. The  desired result then follows. \qed

In  \corref{conv-to-t*},  we used the tightness of the set $\{\th_n\}$, which can be proved using the
argument of \lemref{lem41}.
The result indicates that as the stepsize $\e\to 0$ and $n\to \infty$ with $n\e\to \infty$,
$\th_n$ converges to $\th^*$ in the sense in probability.
Note that if $\th^*$ turns out to the optimum of the search space,
then $\th_n$ converges to the optimum.

\begin{rem}{\rm
Note that for notational simplicity, we have assumed that there is a
unique stationary point of \eqref{ODE}. As far as the convergence is
concerned, one need not assume that there is only one $\theta^*$.
See how multimodal cases can be handled in the related  stochastic approximation problems in \cite[Chpaters 5, 6, 8]{kushner2003}.
In fact, for the multimodal cases, we can show that $\th^\e(\cdot + t_\e)$
converges in an appropriate sense to the set of the stationary
points.
Thus \corref{conv-to-t*} can be modified. In the rate of
convergence study, \cite{Kan89} suggested  an approach using conditional
distribution, which is a modification of a single stationary point.
If multiple stationary points are involved, we can simply use the
approach of \cite{Kan89} combined with our weak
convergence analysis. The notation will be a bit more complex, but main
idea still rest upon the basic analysis method to be presented in
the next section. It seems to be more instructive to present the
main ideas, so we choose the current setting. }\end{rem}

\section{Rate of Convergence}\label{sec:rate}

Once the convergence of a stochastic approximation algorithm is
established, the next task is to ascertain the convergence rate.
To
study the
convergence rate, we take a suitably scaled sequence
$z_n = (\th_n-\th^*)/\e^\al$,
for some $\al>0$. The idea is to choose $\al$ such that $z_n$
converges (in distribution) to a nontrivial limit. The scaling
factor $\al$ together with the asymptotic covariance of the scaled
sequence gives us the rate of convergence. That is, the scaling
tells us the dependence of the estimation error $\th_n-\th^*$ on the
step size, and the asymptotic covariance is a mean of assessing
``goodness'' of the approximation. Here the factor $\al=1/2$ is
used. To some extent, this is dictated by the well-known central
limit theorem. For related work on convergence rate of various
stochastic approximation algorithms, see
\cite{lecuyer1998217,yin199999}.

As mentioned above, by using the definition of the rate of
convergence, we are effectively dealing with convergence in the
distributional sense. In lieu of examining the discrete iteration
directly, we are again taking continuous-time interpolations. Three
assumptions are provided in what follows.

\begin{itemize}
\item[(A4)] The following conditions hold:
\begin{itemize}
\item[(i)] in a neighborhood of $\theta^*$,
$\text{Pr}(\cdot,\eta)$ and $\text{Pg}(\cdot,\eta)$ are
  continuously differentiable
for each $\eta$, and
the second derivatives (w.r.t. $\th$) of
$W(\cdot,r_{1},r_{2},\eta)$ and $P(\cdot,\eta)$
exist and  are continuous.
\item[(ii)] denoting by $E_m$ the conditional expectation on the
$\sigma$-algebra $\F_m=\{\th_0,\,r_{1j},r_{2j}, \eta_j:j< m\}$,
and by $\zeta_\theta$ the first partial derivative w.r.t. $\theta$ of
 $\zeta = W$ or $P$, resp.,
for each positive integer $m$, as $n\to\infty$,
\beq{eq42}
\barray
\ad\frac{1}{n}\sum^{n+m-1}_{j=m}E_m\text{Pr}_\th(\th,\eta_j)\to\lbar{\text{Pr}}_\th(\th) \ \hbox{ in probability,}\\
  \ad\frac{1}{n}\sum^{n+m-1}_{j=m}E_m\text{Pg}_\th(\th,\eta_j)\to\lbar{\text{Pg}}_\th(\th) \ \hbox{ in probability,}\\
 \ad  \sum^\infty_{j=m}|E_m
 W_\th(\th^*,r_{1,j},r_{2,j},\eta_j)|<\infty,\\
\ad \sum^\infty_{j=m}|E_m P_\th(\th^*,\eta_j)-\lbar P_{\theta}
(\theta^*)|<\infty.
\earray
\eeq

\item[(iii)]
The matrix $M+ \lbar P_\theta(\theta^*)$ is stable in that all of
its eigenvalues are on the left half of the complex plane.
\item[(iv)] There is a twice continuously differentiable Lyapunov
function $V\cd: \rr^{2r}\to\rr$ such that
\begin{itemize}
  \item[-] $V(\th)\to\infty$ as $|\th|\to\infty$, and $V_{\th\th}\cd$
  is uniformly
  bounded.
  \item[-] $|V_\th(\th)|\leq K(1+V^{1/2}(\th))$.
  \item[-] $|M\th+\lbar{P}(\th)|^2\leq K(1+V(\th))$ for each $\th$.
  \item[-] $V'_\th(\th)(M\th+\lbar{P}(\th))\leq -\lambda V(\th)$ for some $\lambda>0$ and
  each $\th\neq \th^*$.
 \end{itemize}
\end{itemize}

\item[(A5)]
$\disp\sum^\infty_{j=m} |E\wdt
W'(\theta_m,r_{1,m},r_{2,m},\eta_m)\wdt
W(\theta_j,r_{1,j},r_{2,j},\eta_j)|
\!\!<\infty,$
where $\wdt
W(\theta,r_1,r_2,\eta)=P(\theta,\eta)-\lbar{P}(\theta)+W(\theta,r_1,r_2,\eta).$

\item[(A6)]
The sequence
$B^\e(t)=\sqrt{\e}\sum^{t/\e-1}_{j=0}\wdt W(\th^*,r_{1,j},r_{2,j},\eta_j)$
converges weakly to $B\cd$, a Brownian motion whose covariance
$\Sigma t$ with $\Sigma\in\rr^{2r\times 2r}$ given by
\beq{asy-cov}\barray \Sigma\ad=E\wdt
W(\theta^*,r_{1,0},r_{2,0},\eta_0)\wdt
W'(\theta^*,r_{1,0},r_{2,0},\eta_0)\\
\aad \ \ +\sum^\infty_{k=1}E\wdt
W(\theta^*,r_{1,0},r_{2,0},\eta_0)\wdt
W'(\theta^*,r_{1,k},r_{2,k},\eta_k)\\
\aad\ \ +\sum^\infty_{k=1}E\wdt
W(\theta^*,r_{1,k},r_{2,k},\eta_k)\wdt
W'(\theta^*,r_{1,0},r_{2,0},\eta_0).
\earray\eeq
\end{itemize}

\begin{rem}\label{rem:abouta6}{\rm
Note that (A4)(ii) is another noise condition. The motivation is similar to \remref{rem:about-a2}.
The main difference of  \eqref{eq24-0} and \eqref{eq24} and \eqref{eq42} is that \eqref{eq42} is on the derivative of the functions evaluated at
the point $\theta^*$.
In fact, we only  need
the derivative exists in a neighborhood of this point only. This is because that we are
analyzing the asymptotic normality locally. In view of this condition and condition of $\{r_{i,n}\}$,
\bea \ad\frac{1}{n}\sum^{m+n-1}_{j=m}E_m
W_\th(\th^*,r_{1,j},r_{2,j},\eta_j)\to 0\text{  in probability,}\\
\ad \frac{1}{n}\sum^{m+n-1}_{j=m}E_m P_\th(\th^*,\eta_j)\to
\lbar P_{\theta} (\theta^*) \text{  in probability.}\eea
The traditional PSO algorithms do not allow non-additive noise, here we are treating a more general problem. Nonadditive noise can be allowed.

(A4)(iv) assumes the existence of a Lyapunov function. Only the existence is needed; its precise form
need not be known.
For simplicity, we have assumed the convergence of the scaled sequence to a Brownian motion in (A6);
sufficient conditions are well known; see for example, \cite[Section 7.4]{kushner2003}.
Before proceeding further, we first obtain a moment bound of
$\th_n$.
}\end{rem}

\begin{lem}\label{lem41}
Assume that {\rm(A1)-(A6)} hold. Then there is an $N_\e$ such that
for all $n>N_\e$, $EV(\th_n)=O(\e)$.
\end{lem}

\para{Proof.}
To begin, it can be seen that
\beq{eq44}
\barray
\ad E_nV(\th_{n+1})-V(\th_n)\\[1ex]
\aad\leq -\e\lambda V(\th_n)+\e E_nV'_\th(\th_n)\wdt
W(\th_n,r_{1,n},r_{2,n},\eta_n)\\
\aad\quad+O(\e^2)(1+V(\th_n)),
\earray
\eeq
where $\th^+_n$ is on the line segment joining $\th_n$ and
$\th_{n+1}$. The
 bound in \eqref{eq44} follows from the
growth condition in (A4)(iv), the last inequality follows from (A1).
To proceed, we use the methods of perturbed Lyapunov functions,
which entitles to introduce small perturbations to a Lyapunov
function in order to make desired cancelation. Define a perturbation
\[V^\e_1(\th,n)=\e\sum^\infty_{j=n} E_nV'_\th(\th)\wdt W(\th,r_{1,j},r_{2,j},\eta_j).\]
Note that
\beq{eq45}
|V^\e_1(\th,n)|=K \e(1+V(\th)).
\eeq
Moreover,
\beq{eq46}
\barray
\ad E_nV^\e_1(\th_{n+1},n+1)-V^\e_1(\th_n,n)\\
\aad=O(\e^2)(V(\th_n)+1)-\e E_nV'_\th(\th_n)\wdt W(\th_n,r_{1,n},r_{2,n},\eta_n).
\earray
\eeq
Define
$V^\e(\th,n)=V(\th)+V^\e_1(\th,n).$
Using \eqref{eq44} and \eqref{eq46}, we obtain \beq{pert-L}
E_nV^\e(\th_{n+1},n+1)\leq (1-\e\lambda)
V^\e(\th_n,n)+O(\e^2)(1+V^\e(\th_n,n)).\eeq
Choosing $N_\e$ to be
a positive integer such that
$(1-({\lambda\e}/{2}))^{N_\e}\leq K\e.$
Iterating on the recursion \eqref{pert-L}, taking expectation, and
using the order of magnitude
estimate \eqref{eq45}, we can then obtain
\beq{eq47}
\barray
\ad E V^\e(\th_{n+1},n+1)\\
\aad\leq (1-\e\lambda)E
V^\e(\th_n,n)+O(\e^2)(1+V^\e(\th_n,n))\\
\aad\leq (1-\frac{\e\lambda}{2})^n
EV^\e(\th_0,0)+O(\e)=O(\e).
\earray
\eeq
when $n>N_\e$. The second line of \eqref{eq47} follows from
$1-\lambda\e+O(\e^2)\leq 1-\frac{\lambda\e}{2}$ for sufficiently
small $\e$. Now using \eqref{eq45} again, we also have
$EV(\th_{n+1})=O(\e)$. Thus the desired estimate follows.\qed

Define $z_n=(\th_n-\th^*)/\sqrt{\e}$. Then it is
readily verified that
\beq{eq43}
\barray
z_{n+1}=z_n\ad+\e
(M+\lbar{P}_\th(\th^*))z_n\\
\ad+\sqrt{\e}(P(\th^*,\eta_n)-\lbar{P}(\th^*)+W(\th^*,r_{1,n},r_{2,n},\eta_n))
\\
\ad
+\e(P_\th(\th^*,\eta_n)-\lbar{P}_\th(\th^*)\\
\aad\quad\quad+W_\th(\th^*,r_{1,n},r_{2,n},\eta_n))z_n+o(|z_n|^2).
\earray
\eeq

\begin{cor}\label{cor42}
Assume that {\rm(A1)-(A6)} hold. If the Lyapunov function is locally
quadratic, i.e.,
\[V(\th)=(\th-\th^*)'Q(\th-\th^*)+o(|\th-\th^*|^2).\]
Then $EV(z_n)=O(1)$ for all $n>N_\e$.
\end{cor}

Now we are in a position to study the asymptotic properties through
weak convergence of appropriately interpolated sequence of $z_n$.
Define $z^\e(t)=z_n$ for $t\in[(n-N_\e)\e,(n-N_\e)\e+\e]$.
 We
can introduce a truncation sequence. That is,
in lieu of $z^\e\cd$, we let $N$ be a fixed but otherwise arbitrary large positive integer
and define $z^{\e,N}\cd$ as an $N$-truncation of $z^\e\cd$. That is, it is equal to
$z^\e\cd$
up until the first exit of the process from the sphere $S_N=\{|z|: |z|\le N\}$
with radius $N$. Also define a truncation function $q^N(z)=1$ if $z \in S_N$, $=
0$ if $z \in \rr^r -S_{N+1}$, and is smooth. Corresponding to such a truncation,
we also have a modified operator with truncation (i.e., the functions used in
the operator are all modified by use of $q^N(z)$).
Then we proceed to establish the convergence of $z^{\e,N}\cd$ as a solution of a martingale
problem with the truncated operator.
Then finally, letting $N\to \infty$, we use the uniqueness of the martingale problem to
conclude the proof.
The argument is similar to that of  Section \ref{sec:conv}.
For further technical details, we refer the reader to
\cite[pp. 284-285]{kushner2003}.
 Such a truncation  device is also widely used in the
analysis of partial differential equations. For notational simplicity, we choose
to simply assume the boundedness rather than go with the truncation route.
Thus merely for notational simplicity, we suppose $z^\e(t)$  is
 bounded. For the rate of convergence, our focus is on
the convergence of the sequence $z^\e\cd$. We shall show that it
converges to a diffusion process whose covariance matrix together
with the scaling factor will provide us with the desired convergence
rates. Although more complex than \thmref{conv}, we still use the
martingale problem setup. To keep the presentation relatively brief,
we shall only outline the main steps needed.

For any $t,s>0$,
\beq{eq48}
\barray
z^\e(t+s)-z^\e(t)\ad=\e\sum^{(t+s)/\e-1}_{j=t/\e}(M+\lbar{P}_\th(\th^*))z_j\\
\aad\  +\sqrt{\e}\sum^{(t+s)/\e}_{j=t/\e}\wdt W(\th^*,r_{1,j},r_{2,j},\eta_j)\\
\aad \ +\e\sum^{(t+s)/\e}_{j=t/\e}\wdt
W_\th(\th^*,r_{1,j},r_{2,j},\eta_j)z_j.
\earray
\eeq
Note that for any $\dl>0$, $t,s>0$ with $s<\dl$,
\bea
\ad
E^\e_t\l\sqrt{\e}\sum^{(t+s)/\e-1}_{j=t/\e}\wdt W(\th^*,r_{1,j},r_{2,j},\eta_j)\r^2\\
\aad \ \leq K\e\left(\frac{t+s}{\e}-\frac{t}{\e} \right)=Ks\leq
K\dl.
\eea
and similarly,
\bea
\ad
E^\e_t\l\e\sum^{(t+s)/\e-1}_{j=t/\e}\wdt W_\th(\th^*,r_{1,j},r_{2,j},\eta_j)z_j\r^2
\le Ks\leq K\dl.
\eea
Using \corref{cor42} and similar argument as that of \thmref{conv},
we have the following result.

\begin{lem}
Assume conditions of \corref{cor42}, $\{z^\e\cd\}$ is tight on
$D([0,T]:\rr^{2r})$.
\end{lem}

Next we can extract a convergent subsequence of $\{z^\e\cd\}$.
Without loss of generality, still denote the subsequence by
$z^\e\cd$ with limit $z\cd$.
 For any $t,s>0$, \eqref{eq48} holds.
The way to derive the limit is similar to that of \thmref{conv}
using martingale problem formulation although the analysis is more
involved.
 We proceed to show that the limit is the unique solution for the
 martingale problem with operator
\beq{diff-op} L f(z) = {1\over 2} \hbox{tr}(\Sigma f_{zz}(z)) +
(\nabla f(z))'(M + \lbar P(\theta_*)), \eeq for $ f\in C^2_0$, $C^2$
functions with compact support.

Using similar notation as that of Section \ref{sec:conv}, redefine
\beq{eq311-a}
\wdt h=h(z(t_i):i\leq\kappa),\ \ \wdt h^\e=h(z^{\e}(t_i):
i\leq\kappa).
\eeq
By (A4) (ii), as $\e\to0$
\bea
\ad\!\!\! E \wdt h^\e\bgl \e\sum^{(t+s)/\e}_{j=t/\e}\wdt
W(\th^*,r_{1,j},r_{2,j},\eta_j)z_j\bgr\\
\ad= E\wdt  h^\e \bgl
\e\sum^{(t+s)/\dl_\e}_{l=t/\dl_\e}\sum^{lm_\e+m_\e-1}_{j=lm_\e}\wdt
W(\th^*,r_{1,j},r_{2,j},\eta_j)z_j\bgr
\to 0.
\eea
Using the notation as in Section \ref{sec:conv}, \bea \ad\!\!\! E \wdt
h^\e \bgl
\sum^{(t+s)/\dl_\e}_{l=t/\dl_\e}\!\!\frac{\dl_\e}{m_\e}\!\!\sum^{lm_\e+m_\e-1}_{j=lm_\e}\wdt
W(\th^*,r_{1,j},r_{2,j},\eta_j)[ z_j-z_{lm_\e}]\bgr
\\
\ad  \to 0 \ \hbox{ as } \ \e\to 0.\eea
Moreover, by (A6) we have
\bea \disp \sqrt{\e}\sum^{(t+s)/\e}_{j=t/\e}\wdt
W(\th^*,r_{1,j},r_{2,j},\eta_j)\to \int^{t+s}_t d B(u)
\eea
as $\e\to 0$.
For the first term of \eqref{eq48}, we have
\bea \ad E\wdt h^\e \bgl
\e\sum^{(t+s)/\e-1}_{j=t/\e}(M+\lbar{P}_\th(\th^*)) z_j\bgr\\
\aad \ \to E\wdt
h\bgl\int^{t+s}_{t}(M+\lbar{P}_\th(\th^*))z(u)du
 \bgr\eea
as $\e\to 0$. Putting the aforementioned arguments together,
 we have the following theorem.

\begin{thm}
 Under conditions {\rm(A1)-(A7)}, $\{z^\e\cd\}$ converges to $z\cd$
such
that $z\cd$ is a solution of the following stochastic differential equation
\beq{eq49}
dz=[M+ \lbar P_\theta(\theta^*)]zdt+ \Sigma^{1/2}d\wdh B(t),
\eeq
where $\wdh B\cd$ is a standard Brownian motion.
\end{thm}

\begin{rem}
{{To see what kind of functions and the associated ODE and SDE we
are working with,  we look at two simple examples. In the first
example we use $F(x)=x^2$, take $2$ particles, $\chi=1$,
$\kappa_1=-0.271$, $\kappa_2=1$, $c_1=c_2=1.5$, and assume
$\{\eta_k\}$ is an i.i.d. sequence with mean $[0,0,0,0]'$ and
variance $I$. Then \beq{form-M}
\barray
M
 \ad = \left[ {\begin{array}{*{20}c}
   { - 0.271} & 0 & { - 1.5} & 0  \\
   0 & { - 0.271} & 0 & { - 1.5}  \\
   1 & 0 & { - 1.5} & 0  \\
   0 & 1 & 0 & { - 1.5}  \\
\end{array}} \right],
\earray
\eeq
and the limit ODE is given by
$$\dot \th (t)= M \th(t).$$
Thus
$\theta^*=[0,0,0,0]'$ is the minimizer of the swarm,
and $P_\theta (\theta ^* )=
0\in \rr^{4\times 4}$ (a
$4\times4$ matrix with all entries being 0). In the standard
optimization algorithm, one processor is running to approximate the
optimum. Here, we have two particles running simultaneously.
Note that $\th$ has four components. Two of them represent the particles' positions,
and the other two are the particles'  speeds.
 At the
end, both of the particles reach the minimum, representing something that
might be called ``overlapping.'' In addition, eventually the speeds
of both particles reach 0 (or at resting point). As far as the rate
of convergence is concerned, we conclude that $\theta_n -\theta^*$
decades in the order of $\sqrt \e$ (in the sense of convergence in distribution).
Not only is the mean squares error of $(\th_n -\th^*)$  of the order $\e$, but also
the interpolation of the scaled sequence $z_n$
has a limit represented by
a stochastic differential equation
$$dz = M zdt + d\wdh B(t).$$
That is, \eqref{eq49} is satisfied with
$P_\th(\th^*)=0$ and $\Sigma=I$. As
illustrated in \cite{kushner2003}, the scaling factor $\sqrt \e$
together with stationary covariance of the SDE gives us the rate of
convergence. In terms of the swarm, loosely, we have $\th_n -\th^*
\sim N(0, \e \Xi_0)$ [that is, $(\th_n-\th^*)$ is asymptotically normal with
mean $0 \in \rr^4$ and covariance matrix $ \e \Xi_0$],
where $\Xi_0$ is the asymptotic covariance matrix that is the solution
of the Lyapunov equation $M \Xi_0 + \Xi_0 M' = -I$.}

{Likewise, in the second example,  $F(x)=\sin x$ with $x\in[0,1]$.
We  still take $2$ particles, same parameters setting, and assume
$\{\eta_k\}$ is the same i.i.d. sequence as before. Then
$M$ is as in \eqref{form-M}, and
\bea
P_\theta  (\theta ^* ) = \left[ {\begin{array}{*{20}c}
   {0.75} & 0 & {0.75} & 0  \\
   0 & {0.75} & 0 & {0.75}  \\
   {0.75} & 0 & {0.75} & 0  \\
   0 & {0.75} & 0 & {0.75}  \\
\end{array}} \right]\left[ {\begin{array}{*{20}c}
   1 & 0 & 0 & 0  \\
   0 & 1 & 0 & 0  \\
   1 & 0 & 0 & 0  \\
   0 & 1 & 0 & 0  \\
\end{array}} \right].
\eea
It follows that \eqref{eq49} holds with
$$
M+ \lbar P_\theta(\theta^*)= \left[ {\begin{array}{*{20}c}
   {1.229} & 0 & { - 1.5} & 0  \\
   0 & {1.229} & 0 & { - 1.5}  \\
   {2.5} & 0 & { - 1.5} & 0  \\
   0 & {2.5} & 0 & { - 1.5}  \\
\end{array}} \right]
$$
and $\Sigma= I$.
Similar to the previous example, we have that
$\th_n -\th^* $ is asymptotically normal with mean $0$ and covariance
$\e  \wdt\Xi$, where $\wdt \Xi$ is the asymptotic covariance satisfying the
Lyapunov equation
$(M+ \lbar P_\theta(\theta^*)) \wdt \Xi + \wdt \Xi (M+ \lbar P_\theta(\theta^*))'= -I$.
}
}
\end{rem}

\section{Numerical Examples}\label{sec:exm}

We use two simulation examples to
illustrate
the convergence
properties.
Using \eqref{generalform}, we take $\e=0.01$, $\chi=1$,
$\kappa_1=-0.271$, $\kappa_2=1$, $c_1=c_2=1.5$. For simplicity, we
take the additive noise ${\rm Pr}(\th_n,\eta_n)={\rm
Pr}(\th_n)+\eta_n$ and ${\rm Pg}(\th_n,\eta_n)={\rm
Pg}(\th_n)+\eta_n$, where $\eta_n$ is a sequence of i.i.d. random
variables with a standard normal distribution $\mathcal{N}(0,1)$. In
addition, we set the number of swarms to be 5.

\begin{exm} \label{ex51}   Consider the
sphere function: \beq{eq:f1-def}F_1(x)=\sum^D_{i=1}x^2_i,\eeq where
$D$ is the dimension of the variable $x$. Its global optimum is
$(0,0,\ldots,0)'$. First, the dimension of $X$ is set to be 1.
Figures \ref{fig51} and \ref{fig51a} show the state trajectories
(top) and the centered and scaled errors of the first component
$\th_n^1$ (bottom).  The graphs of $\text{Pr}$ (top) and
$\text{Pg}$ (bottom) are also provided.

\begin{figure}
\centering
\includegraphics[width=5cm]{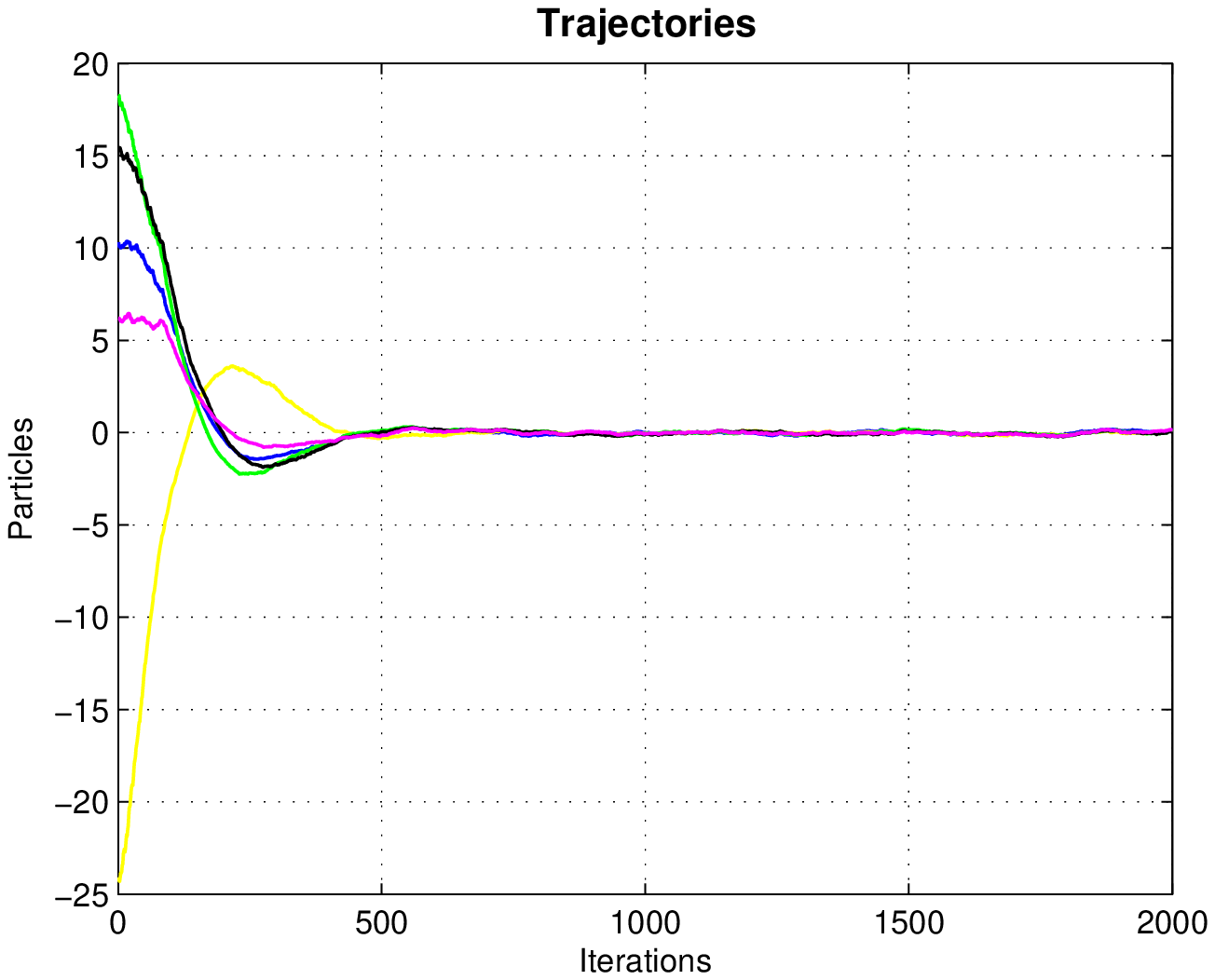}
\includegraphics[width=5cm]{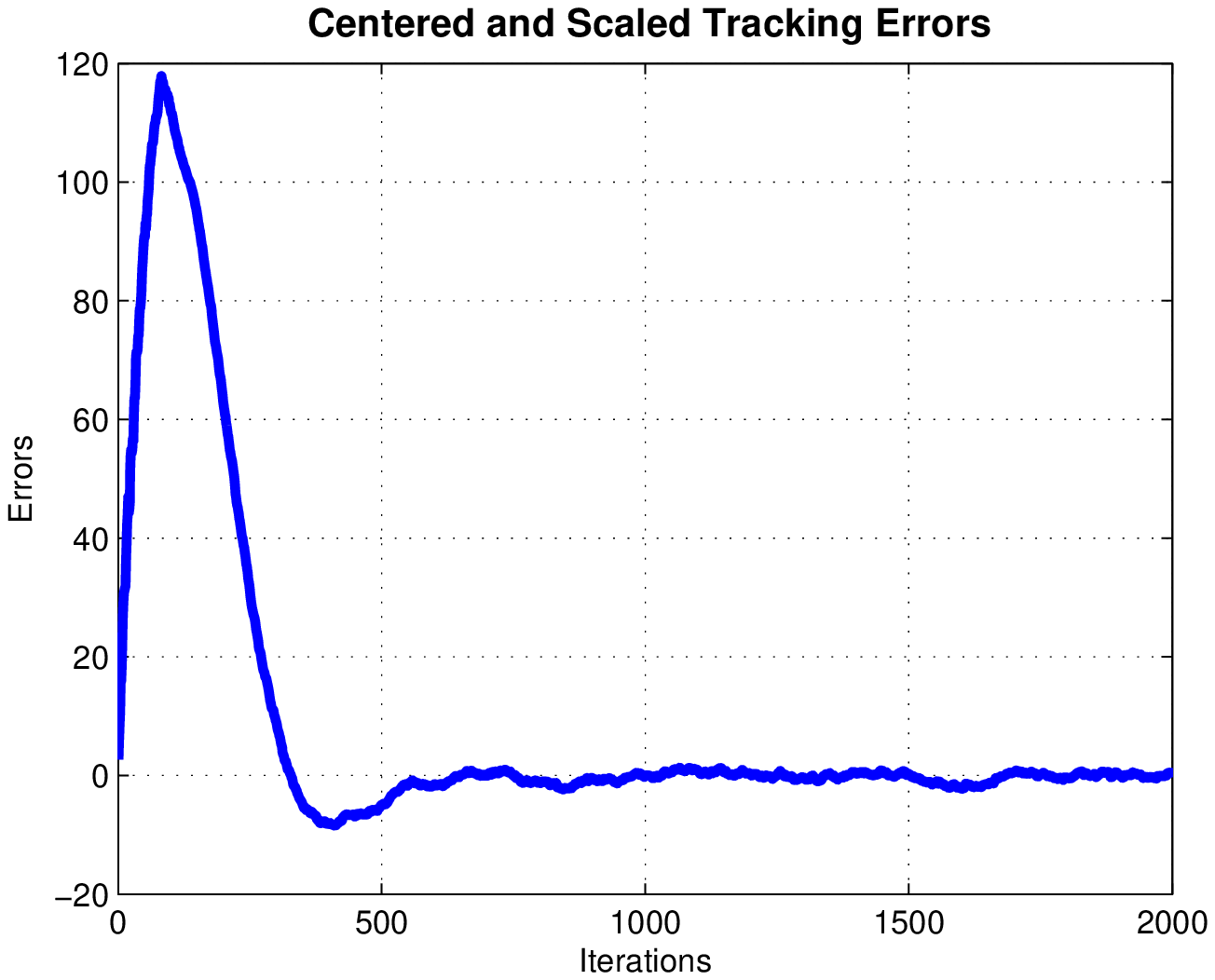}
\caption{Particle swarm of one-dimensional $X$ using $F_1$ defined
in \eqref{eq:f1-def}.} \label{fig51}
\end{figure}

\begin{figure}
  \centering
\includegraphics[width=5cm]{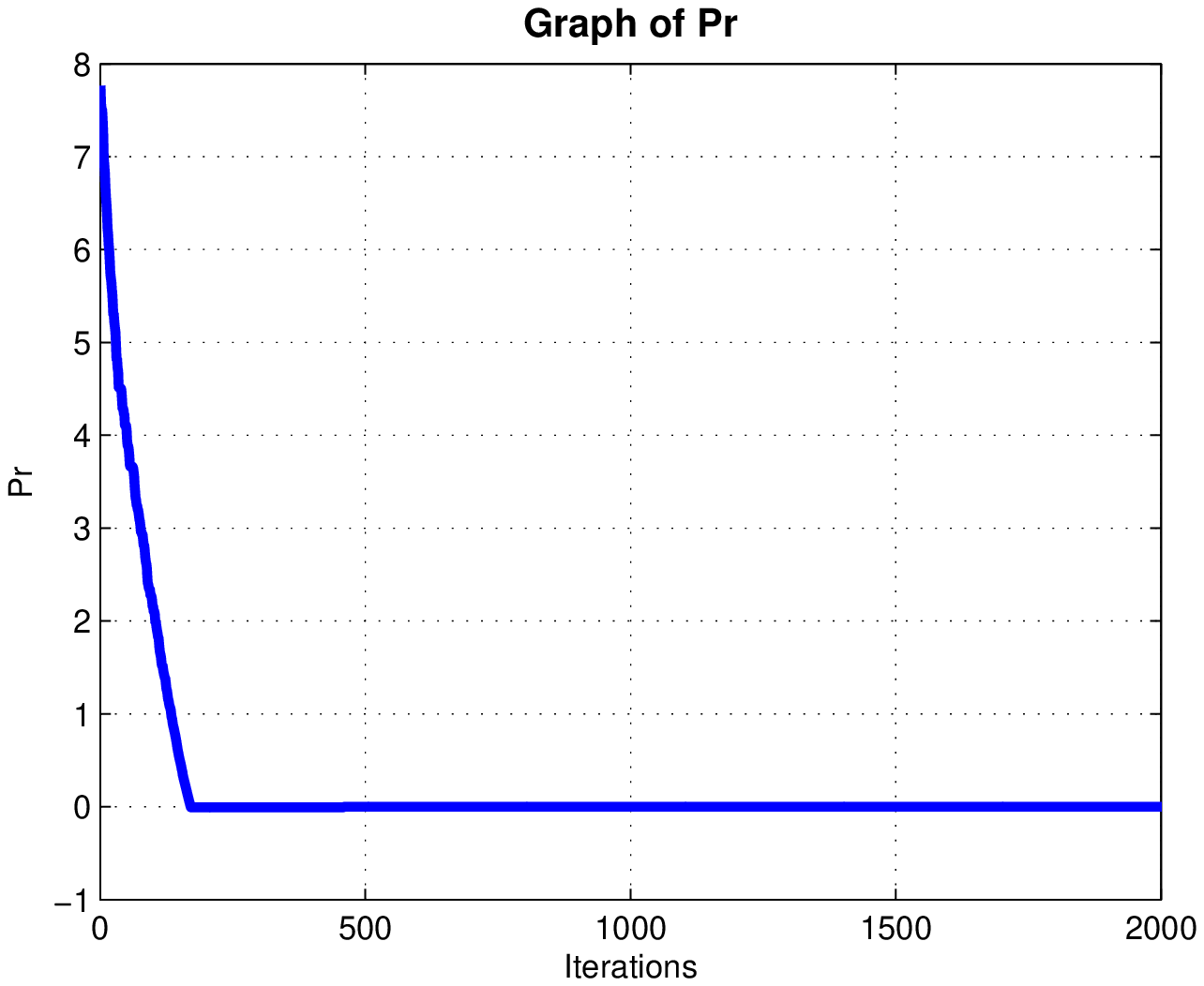}
\includegraphics[width=5cm]{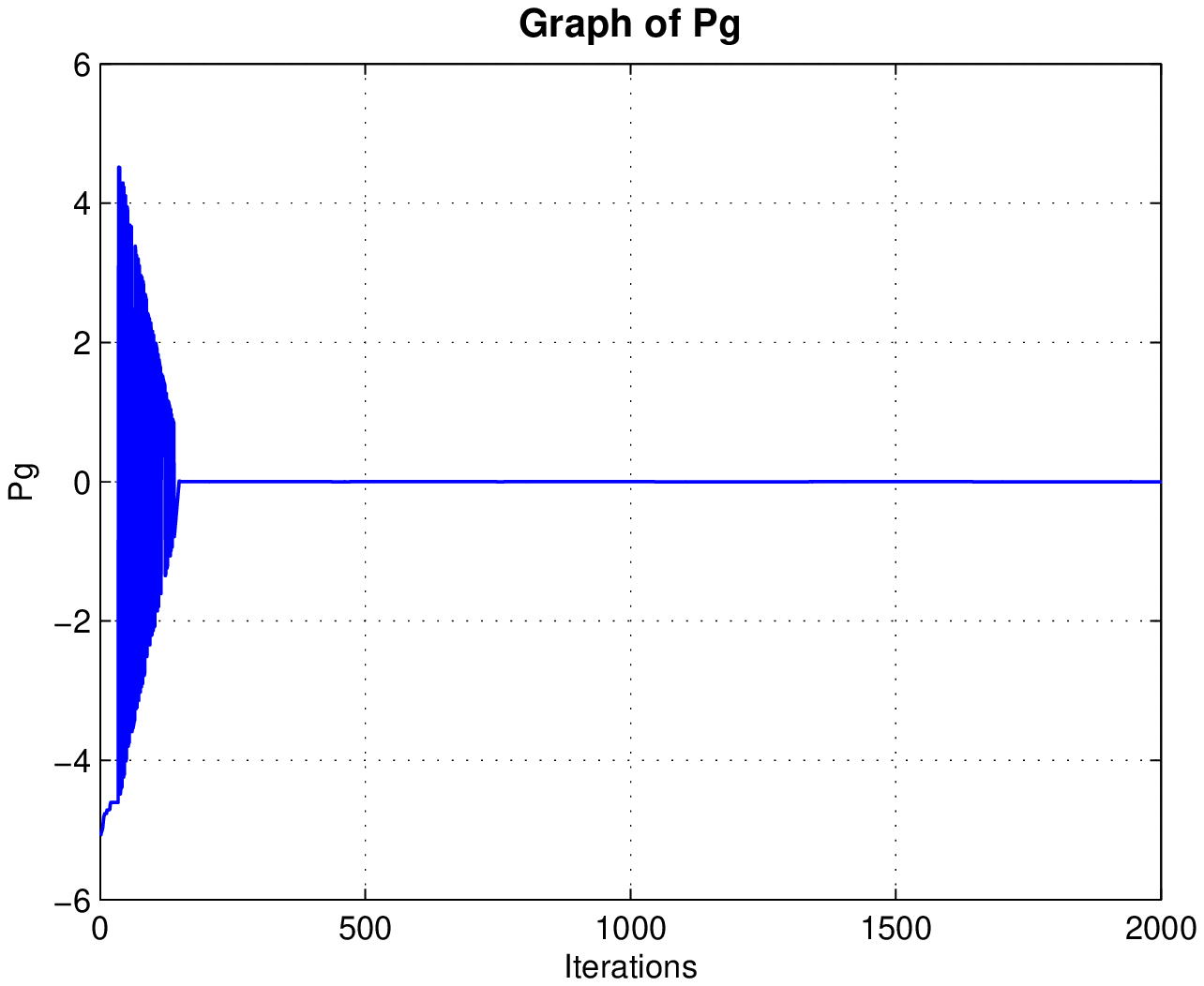}
\caption{Graphs of Pr and Pg using $F_1$ defined in
\eqref{eq:f1-def}.} \label{fig51a}
\end{figure}

Next, we consider the 2-dimension case of $X$. Figures \ref{fig53}
and \ref{fig53a} illustrate the state trajectories (top) and the
centered and scaled errors of the first component $\th_n^1$
(bottom), and the graph of $\text{Pr}$ (top) and $\text{Pg}$
(bottom), respectively.

\begin{figure}
\centering
\includegraphics[width=5cm]{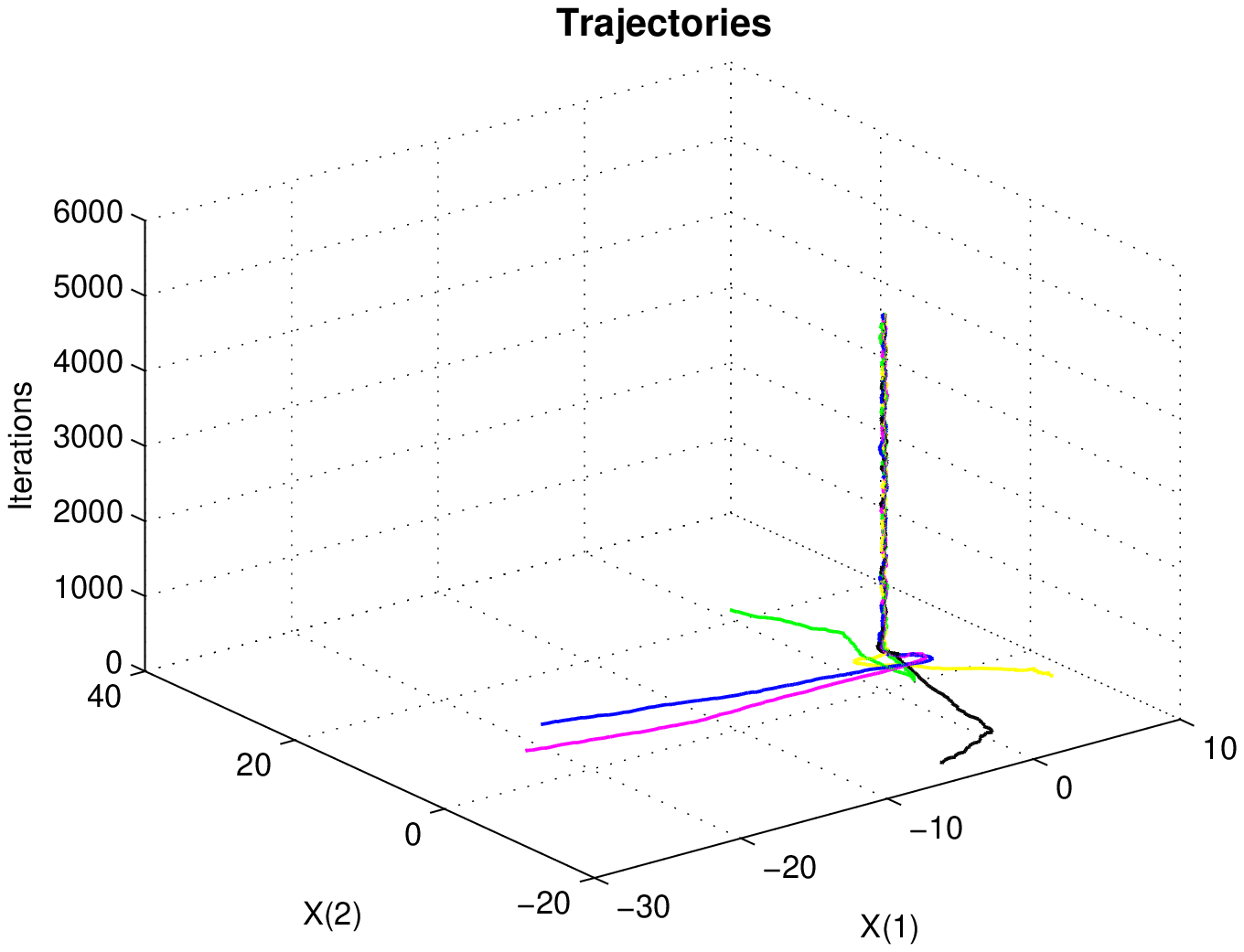}
\includegraphics[width=5cm]{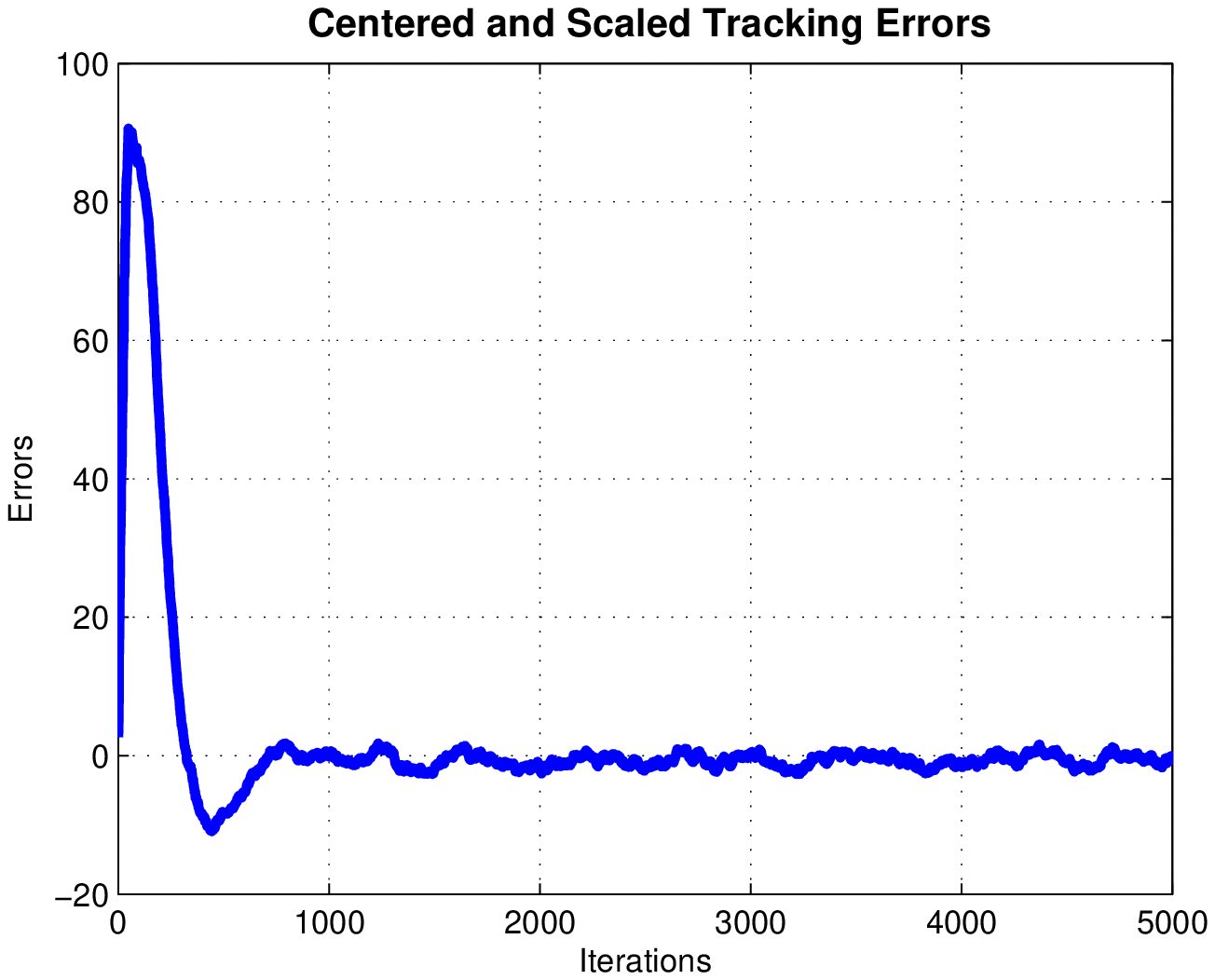}
\caption{Particle swarm of two-dimensional $X$ using $F_1$ defined
in \eqref{eq:f1-def}.} \label{fig53}
\end{figure}
\end{exm}

\begin{figure}
  \centering
\includegraphics[width=5cm]{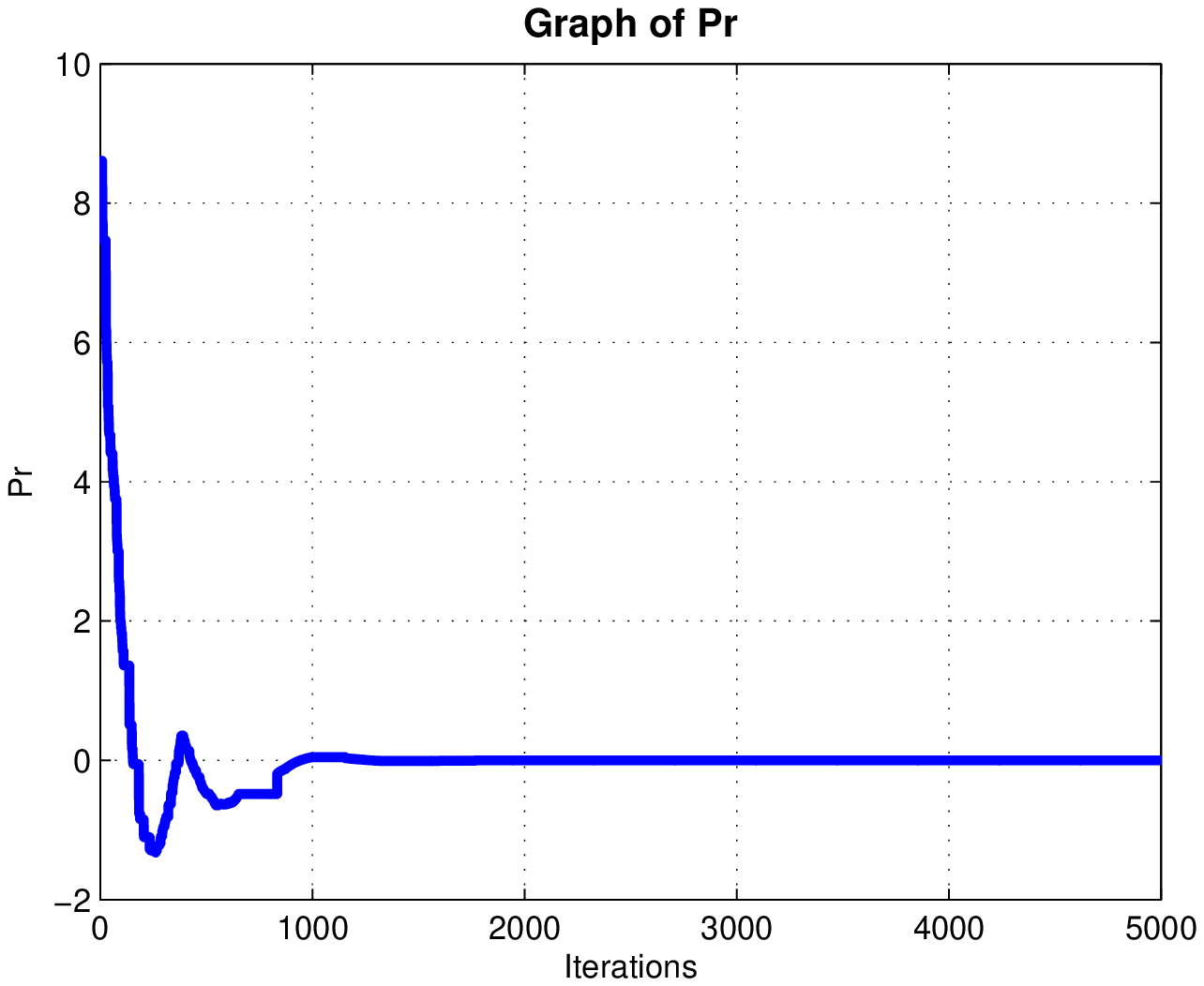}
\includegraphics[width=5cm]{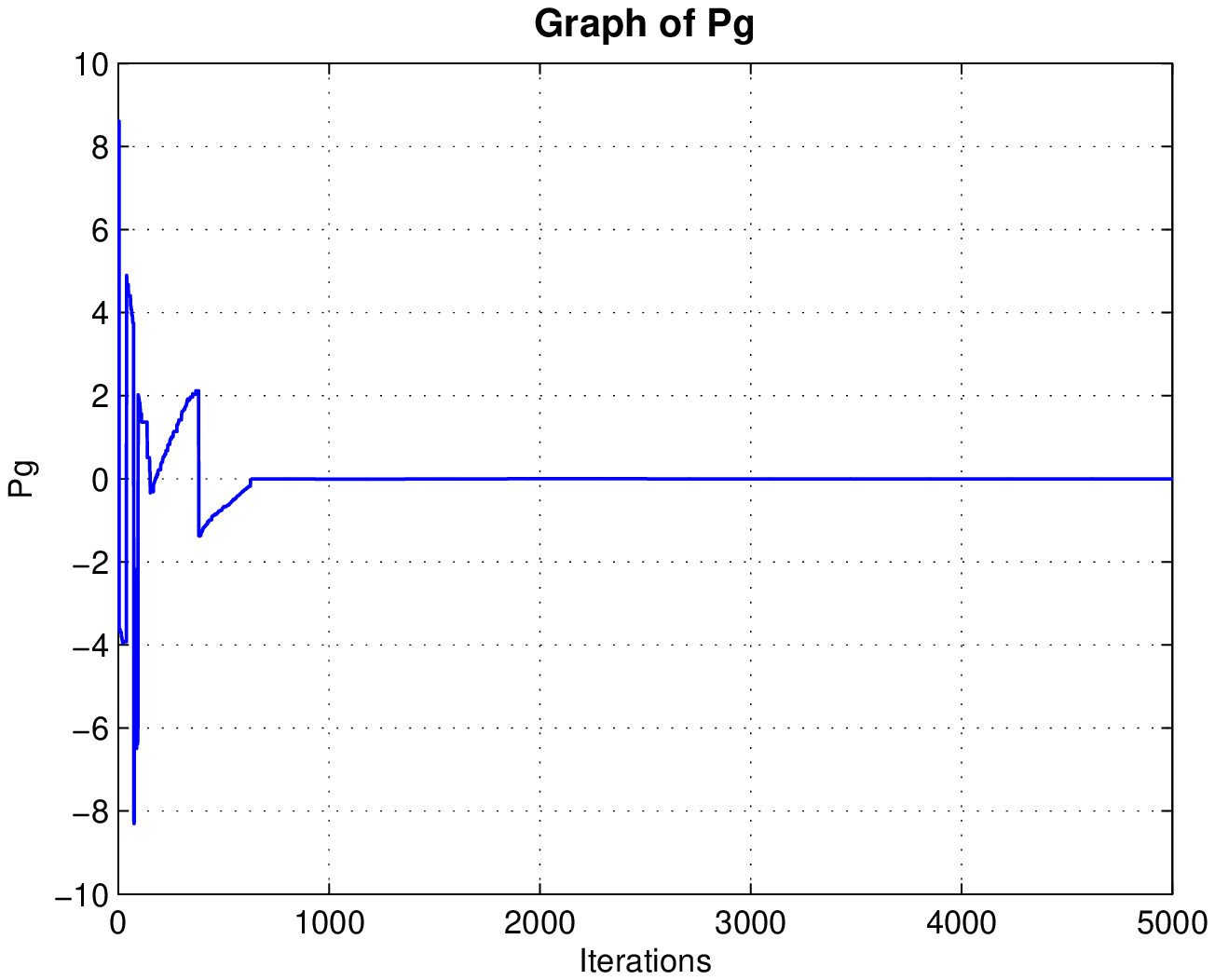}
\caption{Graphs of Pr and Pg using $F_1$ defined in
\eqref{eq:f1-def}.} \label{fig53a}
\end{figure}

\begin{figure}
\centering
\includegraphics[width=5cm]{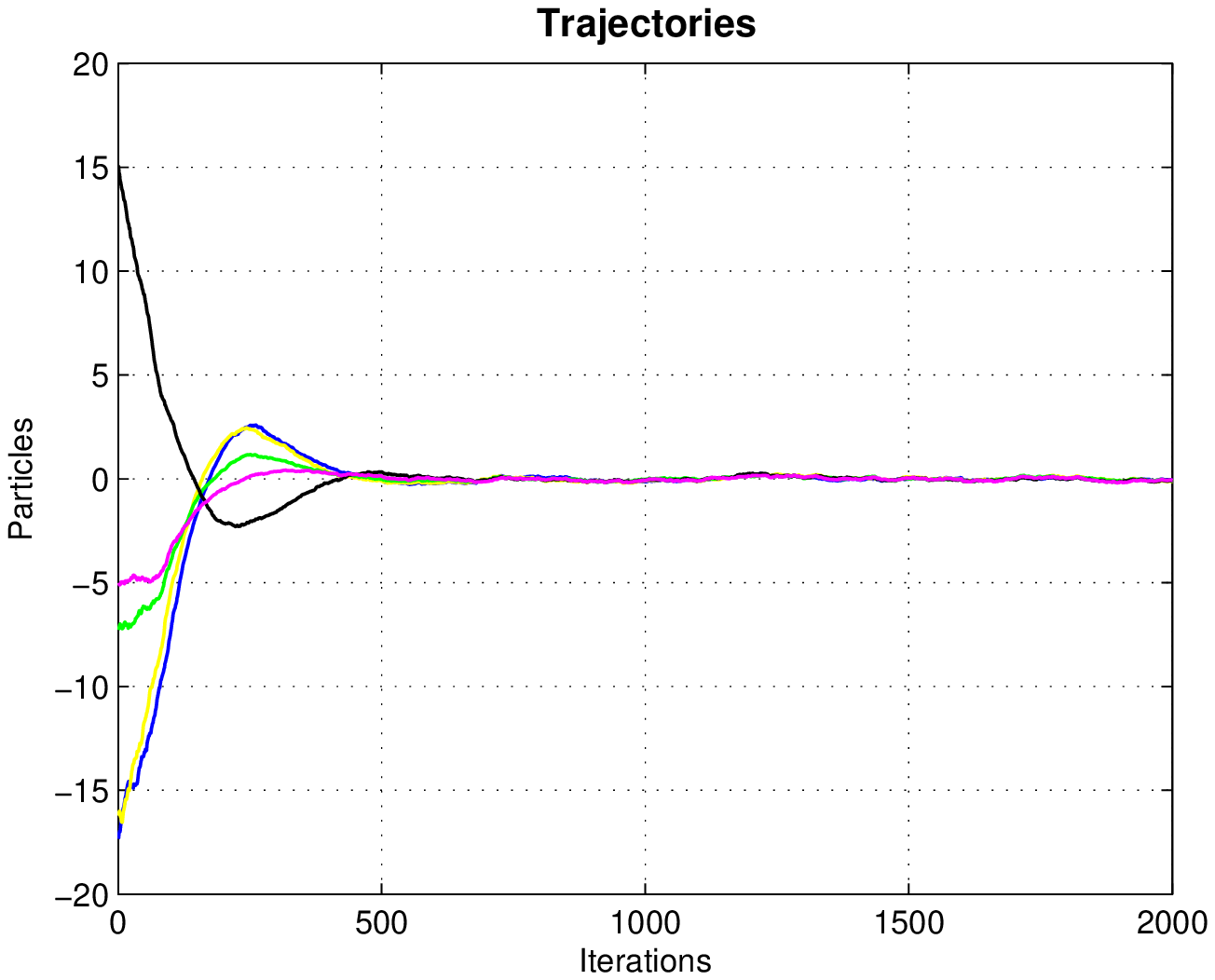}
\includegraphics[width=5cm]{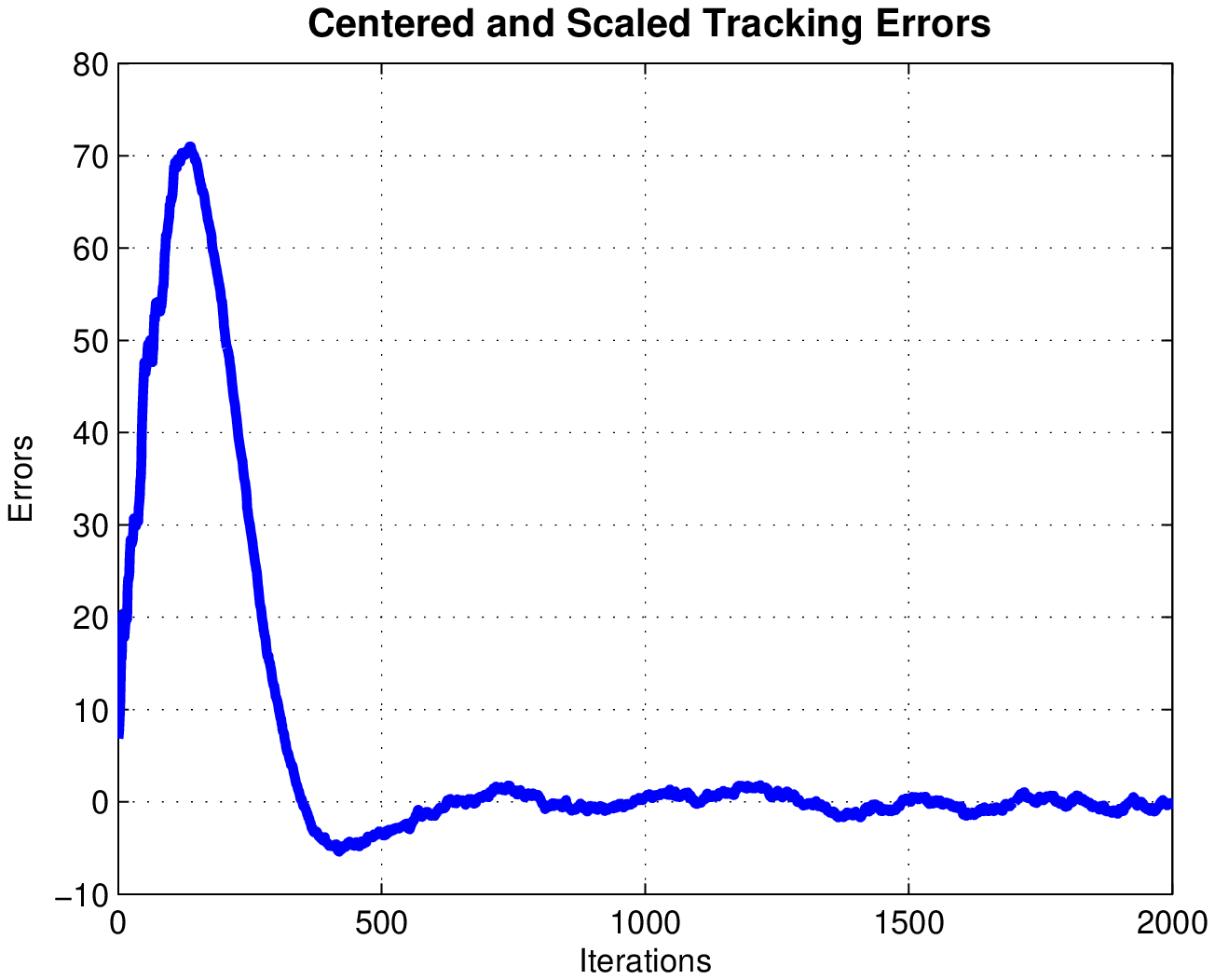}
\caption{Particle swarm of one-dimensional $X$ using $F_2$ defined
in \eqref{eq:f2-def}.} \label{fig52}
\end{figure}

\begin{figure}
  \centering
\includegraphics[width=5cm]{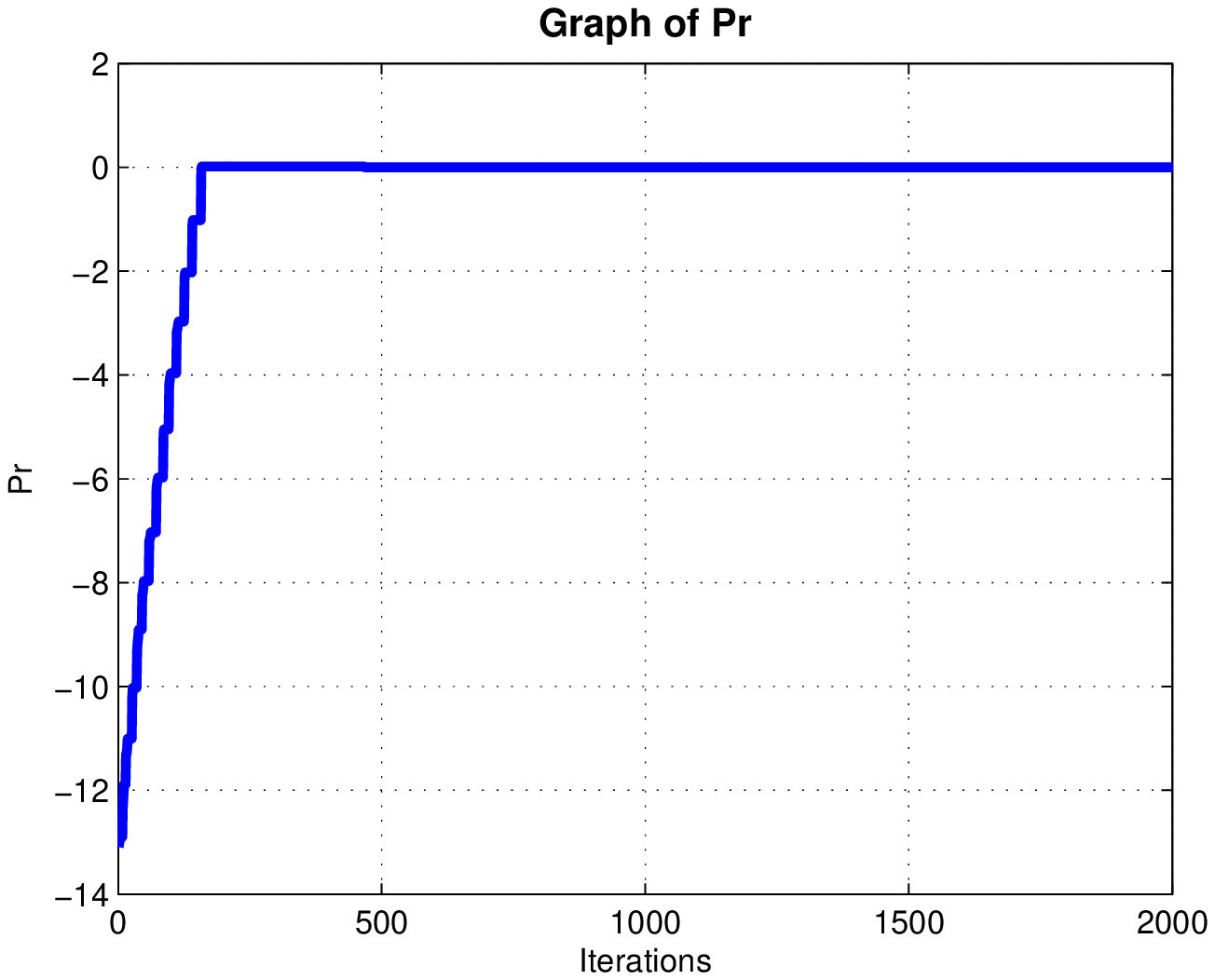}
\includegraphics[width=5cm]{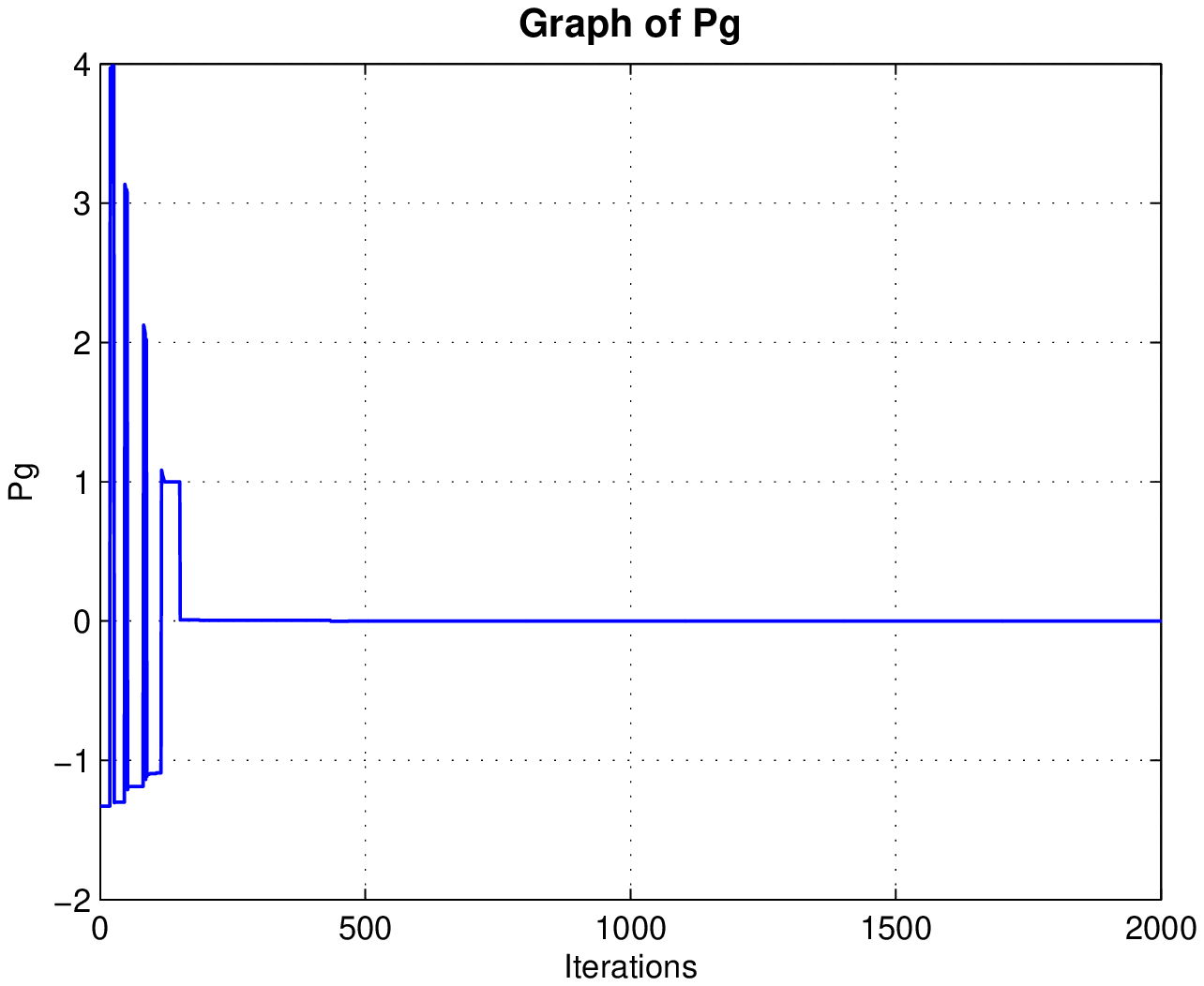}
\caption{Graphs of Pr and Pg using $F_2$ defined in
\eqref{eq:f2-def}.} \label{fig52a}
\end{figure}

\begin{figure}
\centering
\includegraphics[width=5cm]{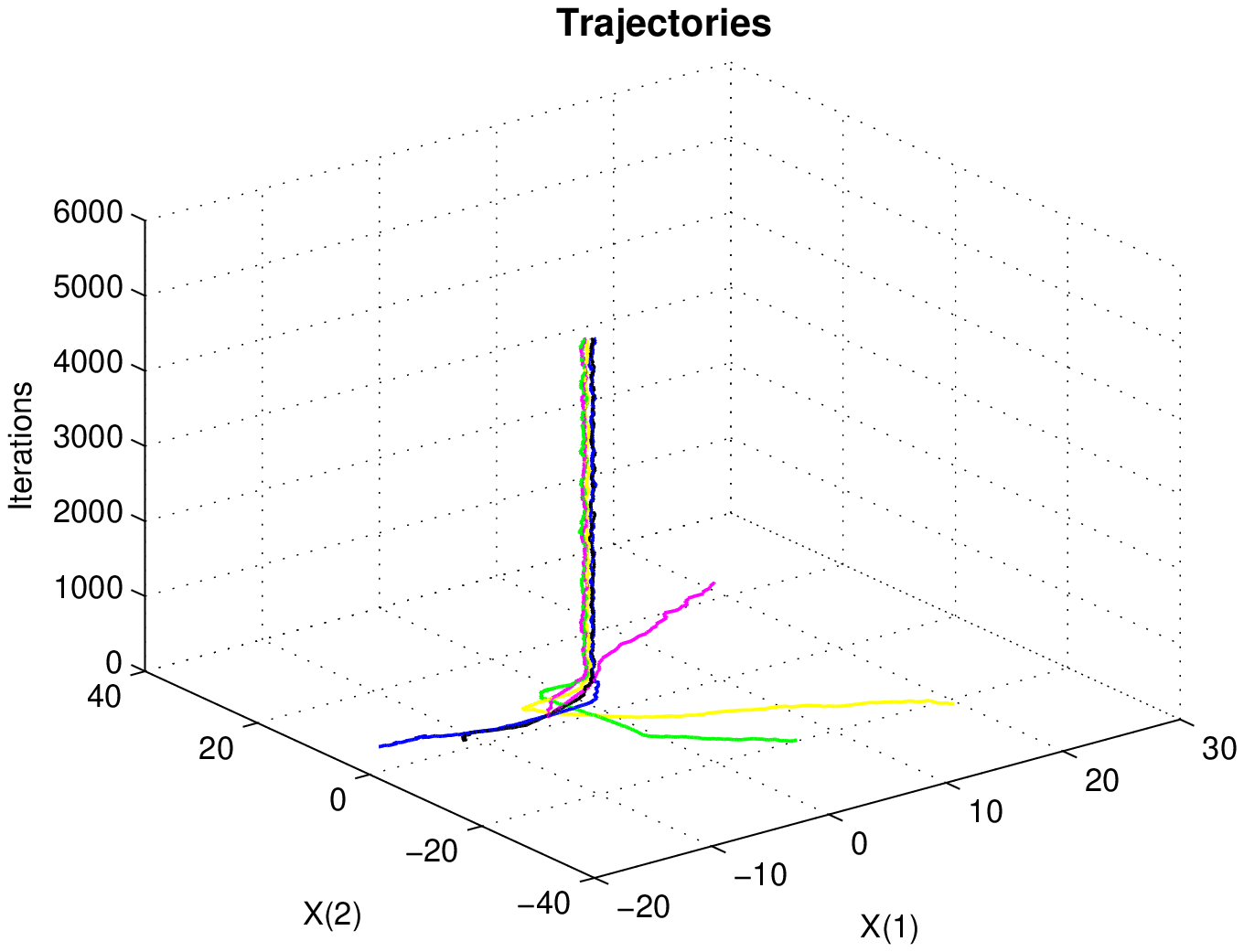}
\includegraphics[width=5cm]{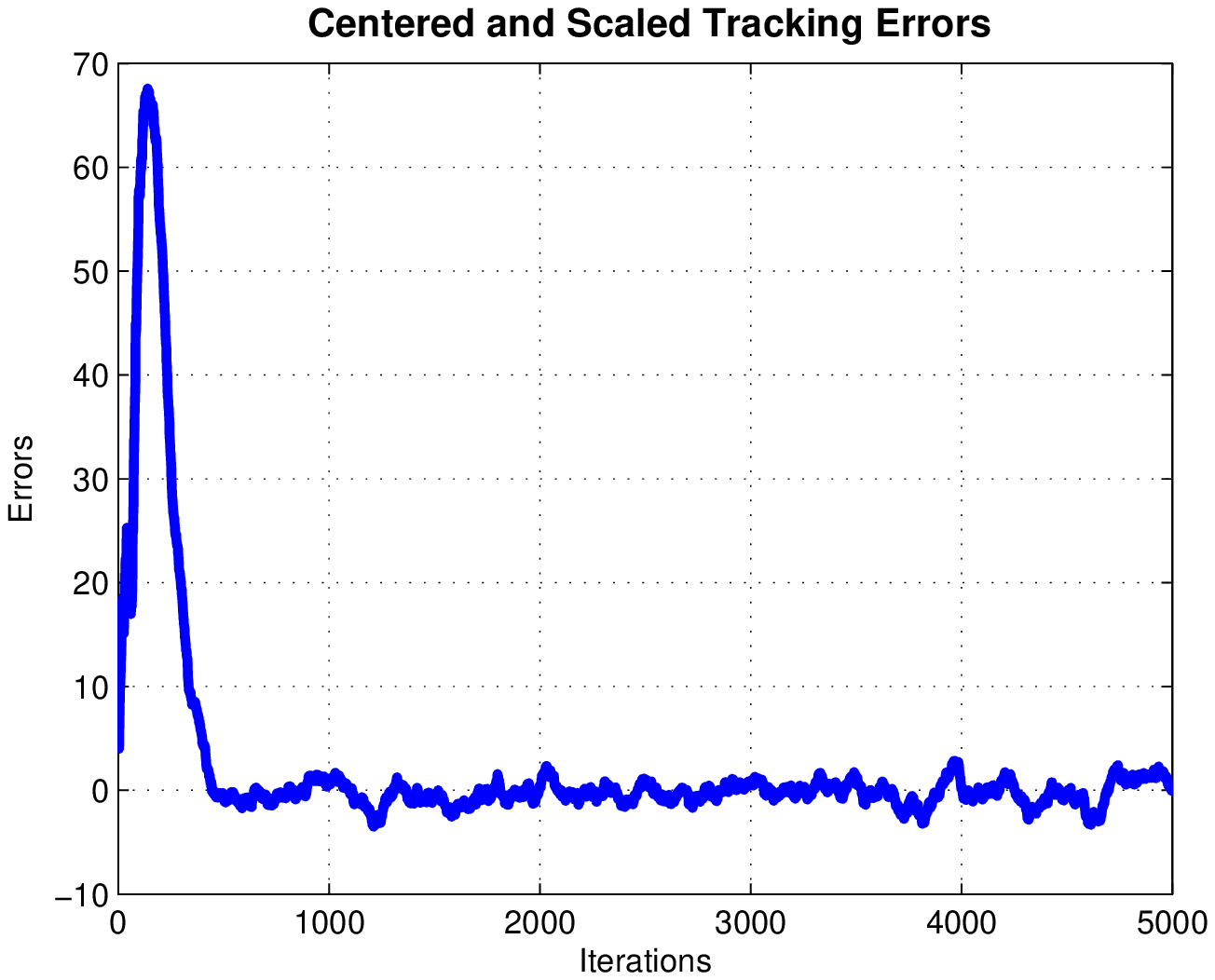}
\caption{Particle swarm of two-dimensional $X$ using $F_2$ defined
in \eqref{eq:f2-def}.} \label{fig54}
\end{figure}

\begin{figure}
  \centering
\includegraphics[width=5cm]{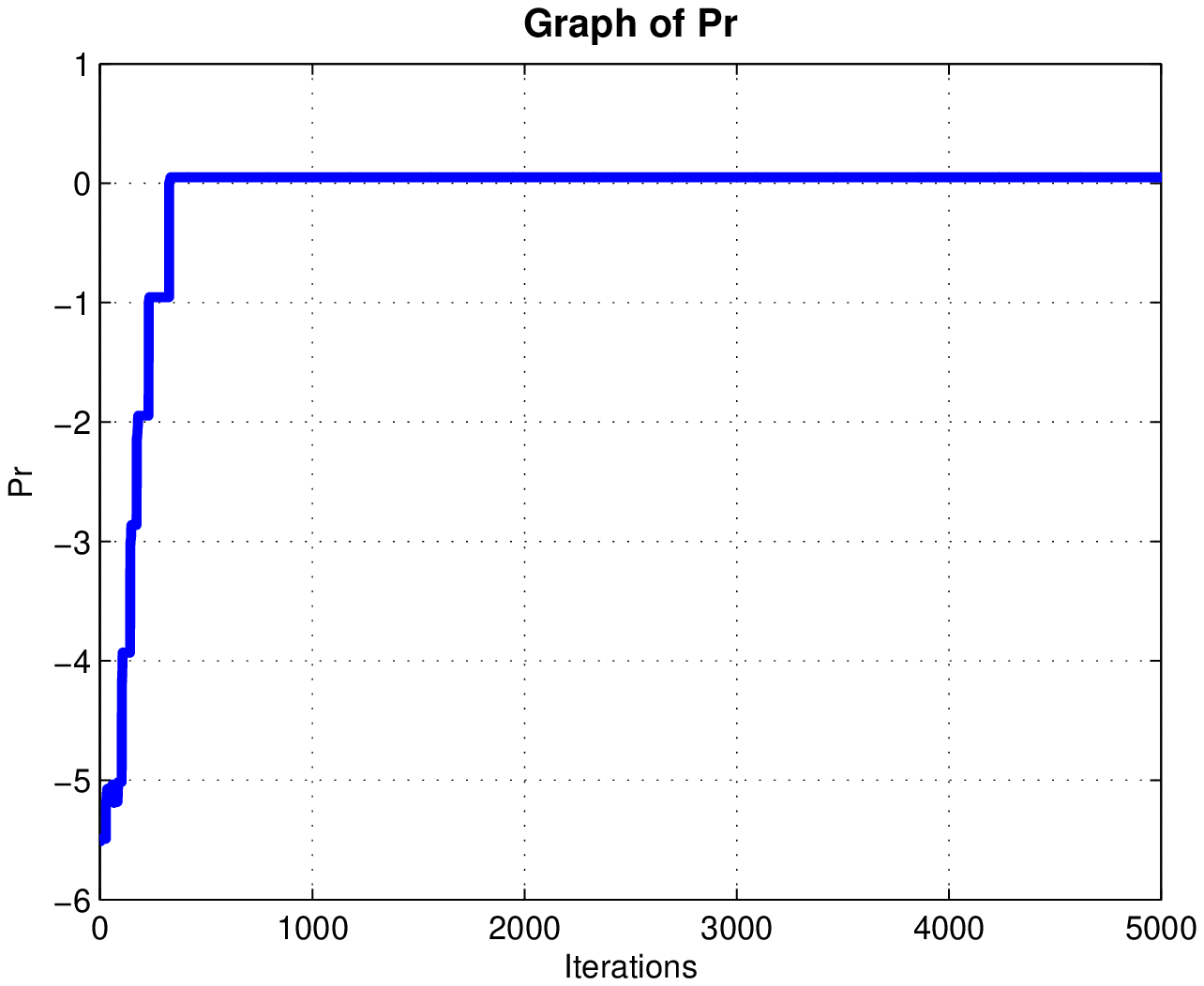}
\includegraphics[width=5cm]{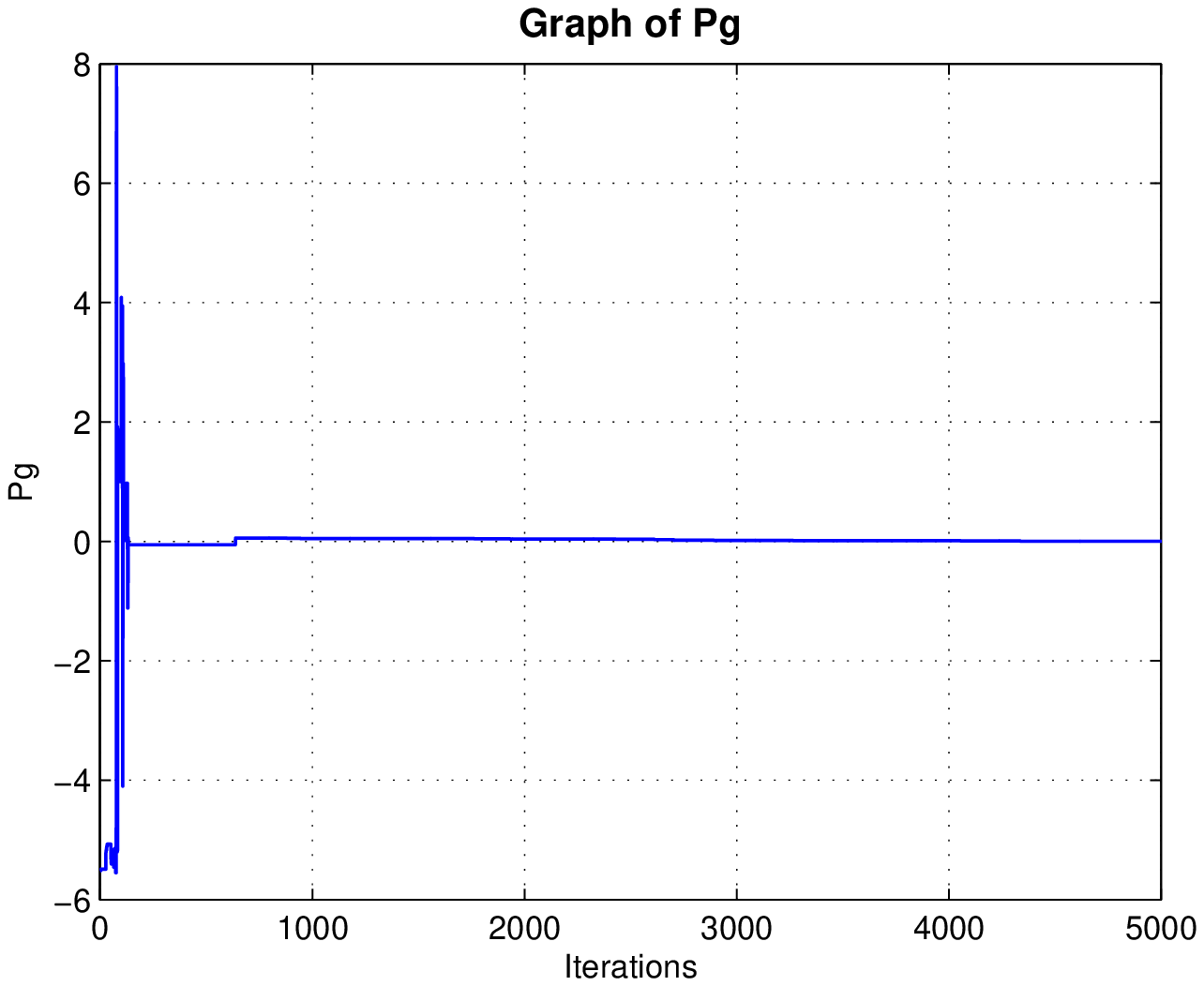}
\caption{Graphs of Pr and Pg using $F_2$ defined in
\eqref{eq:f2-def}.} \label{fig54a}
\end{figure}

\begin{exm} Consider the
Rastrigin function \cite{reynolds199771}
\beq{eq:f2-def}F_2(x)=10D+\sum^D_{i=1}[x^2_i-10\cos(2\pi x_i)],\eeq
where $D$ is the dimension of the variable $x$.
\end{exm}
This function has many local minima. Its global optimum is given by
$(0,0,\ldots,0)'$. Same as \exmref{ex51}, we set the dimension of
$X$ to be 1 and 2, respectively. The particle swarm
trajectories, the centered and scaled errors of the first component,
and graphs of $\text{Pr}$ and $\text{Pg}$ are given in
Figures \ref{fig52} to \ref{fig54a}, respectively.

From these figures, we can conclude that all the swarms converge to
a point in the searching space. These results were obtained without
assuming that
 $r_{1}$, $r_{2}$, $\rm Pr$, and $\rm Pg$ are fixed.
 Our numerical results confirm our
 theoretical
findings in Sections \ref{sec:conv} and \ref{sec:rate}.

\begin{rem}
We use the definition of convergence here that a swarm collapse in
which all particles have converged to a point in the search space.
Sometimes we observe (e.g., in the second example) that the
convergence point is not the global or even local optimum.
This problem, referred to as premature in literatures, occurs
commonly in evolutionary algorithms such as PSOs, genetic
algorithms, evolutionary strategies, etc. Based on our numerical
experiments, we found that if the cost function is unimodal
and with low dimensions,
 the equilibrium coincides with the proper
parameter choice.
The problem of under what conditions the equilibrium coincides with the
optimum deserves to be carefully studied in the future.
\end{rem}

\section{Further Remarks}\label{sec:rem}

{In this paper, we considered a general form of PSO algorithms using
a stochastic approximation scheme. Different from the existing
results in the literature, we have used weaker assumptions and
obtained more general  results without depending on empirical work.
In addition, we obtained rates of convergence for the PSO algorithms
for the first time.}

Several research directions may be pursued in the future. We can
use stochastic approximation methods to analyze other schemes of
PSO, for example, the SPSO2011 considered in
\cite{zambrano20132337}. We can set up a stochastic approximation
similar to \eqref{saform} and analyze its convergence and
convergence rate. Finding ways to systematically choose the
parameter values $\kappa_1$, $\kappa_2$, $c_1$, and $c_2$ is a
practically challenging problem. One thought is to construct a level
two (stochastic) optimization algorithm to select best parameter
value in a suitable sense. To proceed in this direction requires
careful thoughts and consideration.  In addition, we can consider
that some parameters such as $\chi$, $\kappa_1$, etc. are not fixed
but change randomly during iterations or change owing to some random
environment change (for example, see \cite{YinZ10}). The problem to
study is to analyze the convergence and convergence rates in such a
case. Furthermore, using another definition of convergence, i.e.,
the swarm's best known position $\rm Pg$ approaching (converging to)
the optimum of the problem, is another possible study direction.

To conclude, this paper demonstrated convergence properties of a
class of general PSO algorithms and  derived the rates of
convergence
 by using a centered and scaled sequence of  the iterates. This
study opens new arenas for subsequent studies on determining
convergence capabilities of different PSO algorithms and parameters.

\end{document}